\documentclass[a4,11pt]{article}
\usepackage{epsfig}
\usepackage{amsmath}
\usepackage{amssymb}

\newtheorem{theo}{Theorem}[section]

\newtheorem{lema}[theo]{Lemma}
\newtheorem{coro}{Corollary}[theo]

\def\qed{\rule{1.0ex}{1.0ex} \medskip \medskip}

\def \dsp {\displaystyle}

\def\Xint#1{\mathchoice
   {\XXint\displaystyle\textstyle{#1}}%
   {\XXint\textstyle\scriptstyle{#1}}%
   {\XXint\scriptstyle\scriptscriptstyle{#1}}%
   {\XXint\scriptscriptstyle\scriptscriptstyle{#1}}%
   \!\int}
\def\XXint#1#2#3{{\setbox0=\hbox{$#1{#2#3}{\int}$}
     \vcenter{\hbox{$#2#3$}}\kern-.5\wd0}}

\def\dashint{\Xint-}

\def\downbar#1{
\setbox10=\hbox{$#1$}
            \dimen10=\ht10 \advance\dimen10 by 2.5pt
            \ifdim \dimen10<15pt 
               \advance\dimen10 by -0.5pt
               \dimen11=\dimen10
               \advance\dimen10 by 2.5pt
               \lower \dimen11
            \else \lower \ht10 \fi
            \hbox {\hskip 1.5pt \vrule height \dimen10 depth \dp10}}
\def\upbar#1{
\setbox10=\hbox{$#1$}
            \dimen10=\ht10 \advance\dimen10 by \dp10 \advance\dimen10 by 2.5pt
            \ifdim \dimen10<15pt 
               \advance\dimen10 by 2pt \fi
            \raise 2.5pt \hbox {\hskip -1.5pt \vrule height \dimen10}}

\begin{document}

\title{\Large{Para-orthogonal polynomials on the unit circle satisfying three term recurrence formulas}\thanks{This work was support by funds from CAPES, CNPq and FAPESP of Brazil.}}
\author
{
 {C.F. Bracciali$^{a}$, A. Sri Ranga$^{a}$\thanks{ranga@ibilce.unesp.br (corresponding author)},  A. Swaminathan$^{b}$   \,  } \\[1ex]
  {\small $^{a}$Departamento de Matem\'{a}tica Aplicada, IBILCE, }
   {\small  UNESP - Univ. Estadual Paulista} \\
   {\small 15054-000, São José do Rio Preto, SP, Brazil }\\[1ex]
  {\small $^{a}$Departament of Mathematics, }
   {\small  IIT Roorkee } \\
   {\small Roorkee, India }\\[1ex]
}

\date{ }

\maketitle

\thispagestyle{empty}

\begin{abstract}
When a nontrivial measure $\mu$ on the unit circle satisfies the symmetry $d\mu(e^{i(2\pi-\theta)}) = - d\mu(e^{i\theta})$ then the associated OPUC, say $S_n$, are all real.  In this case, in \cite{DelsarteGenin-1986}, Delsarte and Genin have shown that the two sequences  of para-orthogonal polynomials $\{zS_{n}(z) + S_{n}^{\ast}(z)\}$ and $\{zS_{n}(z) - S_{n}^{\ast}(z)\}$ satisfy three term recurrence formulas and have also explored some further consequences of these sequences  of polynomials such as their connections to sequences of orthogonal polynomials on the interval $[-1,1]$.  The same authors, in \cite{DelsarteGenin-1988}, have also provided a means to extend these results to cover any  nontrivial measure on the unit circle. However, only recently in \cite{Costa-Felix-Ranga-2013} and then \cite{Castillo-Costa-Ranga-Veronese-2013} the extension associated with  the para-orthogonal polynomials $zS_{n}(z) - S_{n}^{\ast}(z)$ was thoroughly explored, especially from the point of view of the three term recurrence, and chain sequences play an important part in this exploration.  The main objective of the present manuscript is to provide the theory surrounding the extension associated with the para-orthogonal polynomials $zS_{n}(z) + S_{n}^{\ast}(z)$ for any nontrivial measure on the unit circle. Like in \cite{Costa-Felix-Ranga-2013} and \cite{Castillo-Costa-Ranga-Veronese-2013} chain sequences also play an important role in this theory.  Examples and applications are also provided to justify the results obtained.
\end{abstract}

{\noindent}Keywords: Orthogonal polynomials on the unit circle, Para-orthogonal polynomials, Chain sequences. \\

{\noindent}2010 Mathematics Subject Classification: 42C05, 33C47.


\setcounter{equation}{0}
\section{Introduction} \label{Sec-Intro}

Given a nontrivial probability  measure $\mu(\zeta) = \mu(e^{i\theta})$ on the unit circle $\mathcal{C} = \{\zeta=e^{i\theta}\!\!: \, 0 \leq \theta \leq 2\pi \}$,  the associated sequence $\{S_{n}(z)\}$ of monic orthogonal polynomials on the unit circle (OPUC) are those polynomials with the property
\[
\begin{array}l
\dsp \int_{\mathcal{C}} \overline{S_{m}(\zeta)} S_{n}(\zeta) d\mu(\zeta) = \int_{0}^{2\pi} \overline{S_{m}(e^{i\theta})} S_{n}(e^{i\theta}) d\mu(e^{i\theta}) = \delta_{mn}\kappa_n^{-2}.
\end{array}
\]
The orthonormal polynomials on the unit circle are $s_{n}(z) = \kappa_{n} S_{n}(z)$, $n \geq 0$.  The monic  OPUC satisfy the recurrence
\begin{equation} \label{Szego-A-RR}
\begin{array}l
  S_{n}(z) =  z S_{n-1}(z) - \overline{\alpha}_{n-1}\, S_{n-1}^{\ast}(z), \\[1.5ex]
  S_{n}(z) = (1 - |\alpha_{n-1}|^2) z S_{n-1}(z) - \overline{\alpha}_{n-1} S_{n}^{\ast}(z),
\end{array}
n \geq 1.
\end{equation}
where $\alpha_{n-1} = - \overline{S_{n}(0)}$ and $S_{n}^{\ast}(z) = z^{n} \overline{S_{n}(1/\bar{z})}$.   Following Simon \cite{Simon-book-p1} (see also \cite{Simon-2005}) we will refer to the numbers $\alpha_{n}$, $n \geq 0$, as Verblunsky coefficients.  It is well known that these coefficients  are such that $|\alpha_n| < 1$, $n \geq 0$. It is also known that OPUC are completely determined in terms of these coefficients.

For some recent contributions on this topic we refer to  \cite{BreuRyckSimon-2010, Castil-Garza-Marcell-2011, Castillo-Costa-Ranga-Veronese-2013, Costa-Felix-Ranga-2013, Dam-Mug-Yes-2013, Kheif-Golins-Pehers-Yudit-2011, Pehers-2011, Pehers-Volb-Yudit-2011, Simanek-2012,  Ranga-2010, Tsujimoto-Zhedanov-2009} and references there in. Detailed accounts regarding the earlier research on these polynomials can be found, for example, in Szeg\H{o} \cite{Szego-book}, Geronimus \cite{Geronimus-book-1962}, Freud \cite{Freud-book-1971} and Van Assche \cite{VanAssche-Fields-1997}.  However, for recent and more up to date texts on this subject we refer to the two volumes of Simon \cite{Simon-book-p1} and \cite{Simon-book-p2}.  There is also a nice chapter about these polynomials in Ismail \cite{Ismail-book}.

Para-orthogonal polynomials on the unit circle (POPUC)  associated with the OPUC $\{S_n\}$, given by
\[
     zS_{n-1}(z) - \tau_n S_{n-1}^{\ast}(z), \quad n \geq 1,
\]
where $\{\tau_n\}$ is a any sequence of complex numbers such that $|\tau_n| = 1$, are also  important  in the theory of OPUC.

It is well known that the zeros of OPUC are all within the open unit disk. POPUCs are interesting because their zeros are all simple and lie on the unit circle $|z| =1$.  Moreover,  the interpolatory quadrature rules based on the zeros of these polynomials are those quadrature rules on the unit circle which are analogous to the Gaussian quadrature rules on the real line.  These facts were first observed in  Jones, Nj\aa stad and Thron \cite{JoNjTh-1989} and further explored, for example, in \cite{CantMoraVela-2002, DaNjVA-2003, Golinskii-2002} and many other contributions.

When the measure $\mu$ on the unit circle satisfies the symmetry $d\mu(e^{i(2\pi-\theta)}) = - d\mu(e^{i\theta})$, then the OPUC $S_n$ are all real and, consequently, the real Verblunsky coefficients can be given by $\alpha_{n} = - S_{n+1}(0)$, $n \geq 0$.  In this case, the two sets of para-orthogonal polynomials
\begin{equation} \label{DG-Polynomials-1}
    R_n^{(1)}(z) = zS_{n-1}(z) + S_{n-1}^{\ast}(z) \quad \mbox{and} \quad (z-1)R_n^{(2)}(z) = zS_{n}(z) - S_{n}^{\ast}(z), \quad n\geq 1,
\end{equation} 
considered even earlier than \cite{JoNjTh-1989} by Delsarte and Genin \cite{DelsarteGenin-1986},  are important from the point of view of connecting real OPUC to symmetric orthogonal polynomials on the interval $[-1,1]$.  The importance of this connection has been nicely explored by Zhedanov in \cite{Zhedanov-1998}. A further use of this connection to the problem of frequency analysis, see \cite{BrLiSR-2004}. In \cite{DelsarteGenin-1986},  the polynomials $R_{n}^{(1)}$ and $R_{n}^{(2)}$ referred to as the first and second kind {\em singular predictor polynomials}, respectively, are shown to satisfy the  three term recurrence formulas
\[
   \begin{array}{ll}
      R_{n+1}^{(1)}(z) = (z+1)R_{n}^{(1)}(z) - 4 \, d_{n+1}^{(1)} z R_{n-1}^{(1)}(z), &  \\[-0.5ex]
        & n \geq 1, \\[-0.5ex]
      R_{n+1}^{(2)}(z) = (z+1)R_{n}^{(2)}(z) - 4 \, d_{n+1}^{(2)} z R_{n-1}^{(2)}(z), &
   \end{array}
\]
where $R_{0}^{(1)}(z) = R_{0}^{(2)}(z) = 1$,  $R_{1}^{(1)}(z) = R_{1}^{(2)}(z) = z+1$,
\[
    d_{n+1}^{(1)} = \frac{1}{4} (1 - \alpha_{n-2})(1 + \alpha_{n-1}) \ \ \mbox{and} \ \  \quad d_{n+1}^{(2)} = \frac{1}{4} (1 + \alpha_{n-1})(1 - \alpha_{n}), \quad n\geq 1.
\]
Here, one must take $\alpha_{-1} = -1$.  We will refer to $R_{n}^{(1)}$ and $R_{n}^{(2)}$ as  Delsarte and Genin 1 para-orthogonal polynomials (DG1POP) and Delsarte and Genin 2 para-orthogonal polynomials (DG2POP), respectively.

Since the real Verblunsky coefficients $\alpha_{n} = - S_{n+1}(0)$ are such that $-1 < \alpha_n < 1$, $n \geq 0$, one can easily verify that the sequences $\{d_{n+1}^{(1)}\}_{n=1}^{\infty}$ and $\{d_{n+1}^{(2)}\}_{n=1}^{\infty}$ are both positive chain sequences. For more information on positive chain sequences we refer to Chihara \cite{Chihara-book}.

Delsarete and Genin, in \cite{DelsarteGenin-1988}, also have provided  a means to extend these results, especially those associated with the  DG1POP, to include complex Verblunsky coefficients.   In their analysis the following three term recurrence formula 
\[
     \tilde{R}_{n+1}(z) = (\overline{\beta}_{n} z + \beta_{n}) \tilde{R}_{n}(z)  - z  \tilde{R}_{n-1}(z),  \quad n \geq 1,
\]
plays an important role.   In the present manuscript,  the studies are based on the three term recurrence formula 
\[
    R_{n+1}(z) = \big[(1+ic_{n+1})z+(1-ic_{n+1})\big] R_{n}(z) - 4 d_{n+1} z R_{n-1}(z), \quad n \geq 1,
\]
with $R_{0}(z) = 1$ and $R_{1}(z)= (1+ic_{1})z+(1-ic_{1})$, where $\{c_{n}\}_{n=1}^{\infty}$ is a real sequence and $\{d_{n+1}\}_{n=1}^{\infty}$ is a positive chain sequence.    We remark that the three term recurrence formula considered by  Delsarete and Genin \cite{DelsarteGenin-1988} can also be made equivalent to the above three term recurrence formula. 

Only recently, in \cite{Costa-Felix-Ranga-2013} and \cite{Castillo-Costa-Ranga-Veronese-2013}, the extension associated with the DG2POP to accommodate any nontrivial probability measure was thoroughly explored.   As we have presented in the beginning of section \ref{Sec-RecentDevlop} of the present manuscript, a family of OPUC can be found that lead to the same DG2POP.  The resulting para-orthogonal polynomials are  connected to certain real functions on the interval $[-1,1]$. These  real functions which also satisfy certain orthogonality properties can be viewed as an extension to symmetric orthogonal polynomials on the interval $[-1,1]$ (see \cite{BracMcabePerezRanga-2013}).

The principal objective of the present manuscript is to provide a complete picture of the extension  associated with the DG1POP.  We show that given any nontrivial measure on the unit circle one can find a whole family of  DG1POP and we will  give the moment functionals associated with each of these DG1POP.

The manuscript is organized as follows.  In section \ref{Sec-RecentDevlop}, we  briefly present some results given in  \cite{Costa-Felix-Ranga-2013} and then \cite{Castillo-Costa-Ranga-Veronese-2013} which are important for the development of the new results obtained in sections \ref{Sec-Recovering-NewSzego} and \ref{Sec-OPUCtoDG1POP} of the manuscript. These new results also provide a means to detect, with the use of the Verblunsky coefficients, if the associated measure $\mu$ is such that the Integral $\int_\mathcal{C} \zeta(\zeta-1)^{-1} (1-\zeta)^{-1} d \mu(\zeta)$ exists or not.  In section \ref{Sec-FurtherConseq} we give some further  applications. Finally,  in section \ref{Sec-Examples}, examples having explicit formulas are provided to justify the results in the manuscript.

\setcounter{equation}{0}
\section{A brief look at some recent developments} \label{Sec-RecentDevlop}

Some of the results obtained in \cite{Costa-Felix-Ranga-2013} that are relevant for the understanding of the results in the present manuscript can be summarized as follows.  

Let the nontrivial probability measure $\mu(z)$ on the unit circle be such that it has no positive mass (no pure point) at $z=1$. Let the family of nontrivial probability measures on the unit circle $\mu(\delta;z)$, $0 \leq \delta < 1$,  be given by 
\[
      \int_{0}^{2\pi} f(e^{i\theta}) d \mu(\delta;e^{i\theta}) = (1-\delta)\int_{0}^{2\pi} f(e^{i\theta}) d \mu(e^{i\theta}) + \delta f(1).
\]

Let $S_n(\delta; z)$ be the OPUC with respect to the nontrivial probability measure $\mu(\delta;z)$ on the unit circle. Then  for the (appropriately scaled) POPUC  $(z-1)R_n^{(2)}(z)$,  which are independent of the mass $\delta$ at $z=1$, given by
\[
   R_n^{(2)}(z) = \frac{1}{z-1} \frac{\prod_{j=0}^{n-1} \big[1- \rho_{j}^{(2)}\alpha_{j}(\delta)\big]}{\prod_{j=0}^{n-1} \big[1-\mathcal{R}e\big(\rho_{j}^{(2)}\alpha_{j}(\delta)\big)\big]} \Big[ zS_{n}(\delta; z) - \rho_{n}^{(2)}S_{n}^{\ast}(\delta; z)\Big], \quad n \geq 0,
\]
where  $\rho_{0}^{(2)}(\delta) = S_0(\delta;1)/S_0^{\ast}(\delta;1) = 1$ and
\begin{equation} \label{Eq-RhoDG2POP}
       \alpha_{n-1}(\delta) = - \overline{S_{n}(\delta; 0)}  \quad  \rho_{n}^{(2)} = \frac{S_n(\delta; 1)}{S_n^{\ast}(\delta; 1)}
           = \frac{\rho_{n-1}^{(2)} - \overline{\alpha}_{n-1}(\delta)}{1 - \rho_{n-1}^{(2)}  \alpha_{n-1}(\delta)} ,  \quad n \geq 1,
\end{equation}
the following three term recurrence formula holds.
\begin{equation} \label{Eq-TTRR-DG2Ex}
   R_{n+1}^{(2)}(z) = \big[(1+ic_{n+1}^{(2)})z + (1-ic_{n+1}^{(2)})\big] R_{n}^{(2)}(z) - 4\,d_{n+1}^{(2)} z R_{n-1}^{(2)}(z),
\end{equation}
with $R_{0}^{(2)}(z) = 1$ and $R_{1}^{(2)}(z) = (1+ic_{1}^{(2)})z + (1-ic_{1}^{(2)})$, where the real sequences $\{c_n^{(2)}\}$ and the positive chain sequence $\{d_{n+1}^{(2)}\}$ are such that
\begin{equation*}\label{Eq-CoeffsTTRR-2}
  \begin{array}{l}
    \dsp c_{n}^{(2)} = \frac{-\mathcal{I}m \big(\rho_{n-1}^{(2)}\alpha_{n-1}(\delta)\big)} {1-\mathcal{R}e\big(\rho_{n-1}^{(2)}\alpha_{n-1}(\delta)\big)}  = i \,\frac{\rho_{n}^{(2)}-\rho_{n-1}^{(2)}}{\rho_{n}^{(2)}+\rho_{n-1}^{(2)}} \quad \mbox{and} \quad d_{n+1}^{(2)} = d_{1,n}^{(2)} = \big(1-g_{1,n-1}^{(2)}(\delta)\big)g_{1,n}^{(2)}(\delta), 
  \end{array} 
\end{equation*}
for $n \geq 1$. Here, the parameter sequence $\{g_{1,n}^{(2)}(\delta)\}_{n=0}^{\infty}$  of $\{d_{n+1}^{(2)}\}$ is given by
\[
      g_{1,n}^{(2)}(\delta) = \frac{1}{2} \frac{\big|1 - \rho_{n}^{(2)} \alpha_{n}(\delta)\big|^2}{\big[1 - \mathcal{R}e\big(\rho_{n}^{(2)}\alpha_{n}(\delta)\big)\big]}, \quad n \geq 0.
\]
Since $0 < g_{1,0}^{(2)}(\delta)<1$, one can consider the positive chain sequence $\{d_{n}^{(2)}\}_{n=1}^{\infty}$ with the additional element $d_1^{(2)} = d_1^{(2)}(\delta)$ given by $d_1^{(2)} = g_{1,0}^{(2)}(\delta)$.   Then the sequence $\{c_n^{(2)}\}_{n=1}^{\infty}$ together with the minimal parameter sequence $\{m_n^{(2)}\}_{n=0}^{\infty}$ of the positive chain sequence $\{d_{n}^{(2)}\}_{n=1}^{\infty}$, where
\[
    m_0^{(2)} = 0 \quad \mbox{and} \quad m_n^{(2)} = g_{1,n-1}^{(2)}(\delta), \ \  n\geq 1,
\]
can be used to  completely characterize the above measure $\mu(\delta; z)$.  For example, the associated Verblunsky coefficients can be given as
\[
     \alpha_{n-1}(\delta) = \frac{1}{\rho_{n}^{(2)}}\,\frac{1-2m_{n}^{(2)} - i c_n^{(2)}}{1 + i c_n^{(2)}} \quad \mbox{and} \quad \rho_{n}^{(2)} = \frac{1-ic_n^{(2)}}{1+ic_n^{(2)}}\, \rho_{n-1}^{(2)}, \quad n \geq 1,
\]
where $\rho_{0}^{(2)} = 1$. The maximal parameter sequence $\{M_n^{(2)}\}_{n=0}^{\infty}$ of the positive chain sequence $\{d_{n}^{(2)}\}_{n=1}^{\infty}$ is also such that $M_0 = \delta$.

With respect to any of the measures $\mu(\delta;z)$ the polynomials $\{R_n^{(2)}\}$ satisfy the so called L-orthogonality
\[
   \int_{\mathcal{C}} \zeta^{-n+j} R_{n}^{(2)}(\zeta) \,(1-\zeta)d \mu(\delta; \zeta) = 0, \quad 0 \leq j \leq n-1.
\]

The polynomials $R_n^{(2)}(z)$ are actually constant multiples of the CD kernels $K_{n}(\delta; z,1) = \linebreak  \sum_{j=0}^{n} \overline{s_j(\delta;1)}\,s_j(\delta; z)$ associated with the probability measure $\mu(\delta;z)$. Here, $s_n(\delta;z)$ are the orthonormal polynomials associated with $\mu(\delta;z)$.  More on  studies that use the connection between CD kernels  and para-orthogonal polynomials we refer to  Cantero, Moral and Vel\'{a}zquez \cite{CantMoraVela-2002} and Golinskii \cite{Golinskii-2002}. For a much recent paper exploring this connection see \cite{Wong-2007}.  We also cite \cite{Simon-2011}, where there is a nice section on  para-orthogonal polynomials.

As shown in \cite{Castillo-Costa-Ranga-Veronese-2013}, given the real sequence $\{c_n^{(2)}\}_{n=1}^{\infty}$ and the positive chain sequence $\{d_{n}^{(2)}\}_{n=1}^{\infty}$ one can also recover information about the measure $\mu(\delta;z)$ and the associated OPUC $S_n(\delta;z)$ directly with the use of the three term recurrence formula (\ref{Eq-TTRR-DG2Ex}). Some of the results found in \cite{Castillo-Costa-Ranga-Veronese-2013} are briefly stated below. These results will be used in Section \ref{Sec-Recovering-NewSzego} to recover another measure from the same three term recurrence formula. These new results give us the means to answer the question about the extension to the DG1POP.

Let the polynomials $\{Q_n\}$ and $\{R_n\}$ be such that
\begin{equation} \label{Eq-TTRR-Rn}
  \begin{array}l
    Q_{n+1}(z) = \big[(1+ic_{n+1})z+(1-ic_{n+1})\big] Q_{n}(z) - 4 d_{n+1} z Q_{n-1}(z), \\[2ex]
    R_{n+1}(z) = \big[(1+ic_{n+1})z+(1-ic_{n+1})\big] R_{n}(z) - 4 d_{n+1} z R_{n-1}(z),
  \end{array}
    \quad n \geq 1,
\end{equation}
with $Q_0(z) = 0$, $R_{0}(z) = 1$, $Q_1(z) = 2d_1$ and $R_{1}(z)= (1+ic_{1})z+(1-ic_{1})$, where we always assume $\{c_{n}\}_{n=1}^{\infty}$ to be any real sequence.  Any of the coefficients $c_n$ can also be allowed to be zero.  However, when $c_n=0$, $n \geq 1$, then the results will lead to those found in  \cite{BrLiSR-2004}, \cite{DelsarteGenin-1986}  and     \cite{Zhedanov-1998}.

From the recurrence formula (\ref{Eq-TTRR-Rn})
\[
       R_{n}^{\ast}(z) = z^n \overline{R_{n}(1/\bar{z})} = R_{n}(z) \quad \mbox{and} \quad Q_{n}^{\ast}(z) = z^{n-1} \overline{Q_{n}(1/\bar{z})} = Q_{n}(z), \quad n \geq 1.
\]
With these property the polynomials $R_n$ and $Q_n$ can be called self-inversive polynomials or, more appropriately, conjugate reciprocal polynomials. Moreover, if $R_{n}(z) = \sum_{j=0}^{n} r_{n,j}\,z^j$ and $Q_{n}(z) = \sum_{j=0}^{n-1} q_{n,j}\,z^j$, then
\[
   r_{n,n} = \overline{r}_{n,0} = \prod_{k=1}^{n} (1 + i c_k),  \  n \geq 1 \quad \mbox{and} \quad  q_{n,n-1} = \overline{q}_{n,0} = 2d_1 \, \prod_{k=2}^{n} (1 + i c_k), \  n \geq 2.
\]

Firstly, assuming $\{d_{n}\}_{n=1}^{\infty}$ to be any sequences such that $d_n \neq 0$, $n \geq 1$, we give the following lemma that follows from the series expansions of the rational functions $Q_{n}(z)/R_{n}(z)$.

\begin{lema} \label{Lemma-Orto-Rn}
Given the real sequences $\{c_{n}\}_{n=1}^{\infty}$ and $\{d_{n}\}_{n=1}^{\infty}$, where  $d_n \neq 0$, $n \geq 1$, let  the sequences of polynomials  $\{Q_n\}$ and $\{R_n\}$ be as in $(\ref{Eq-TTRR-Rn})$. Then there exist two series expansions
\begin{equation} \label{Eq-Two-Series}
  E_{0}(z) = -\sum_{n=0}^{\infty}\nu_{n+1}\,z^{n}  \quad  \mbox{and} \quad E_{\infty}(z) = \sum_{n=1}^{\infty} \nu_{-n+1}z^{-n} ,
\end{equation}
where
\begin{equation} \label{Eq-ModifiedMoments-Symmetry}
      \nu_{n} = - \overline{\nu}_{-n+1}, \quad n=1,2, \ldots \ ,
\end{equation}
such that there hold the correspondence properties
\begin{equation} \label{origin-infinity-correspondence-1}
    E_{0}(z) - \frac{Q_{n}(z)}{R_{n}(z)} = \frac{\overline{\gamma}_{n}z^{n}}{\overline{r}_{n,n}} + O\big(z^{n+1}\big) \quad \mbox{and} \quad E_{\infty}(z) - \frac{Q_{n}(z)}{R_{n}(z)} = \frac{\gamma_{n}}{r_{n,n}z^{n+1}} + O\big((1/z)^{n+2}\big),
\end{equation}
for $n \geq 0$.  Moreover, if the (moment) functional $\mathcal{N}$ on the space of Laurent polynomials is defined by $\mathcal{N}[z^{-n}] = \nu_{n}$, \ $n=0, \pm1,\pm2, \ldots ,$  then the polynomials $R_n$  satisfy the orthogonality property
\begin{equation*} \label{LOrthogonality-Rn}
    \mathcal{N}[z^{-n+j} R_{n}(z)] = \left\{
      \begin{array}{cl}
        -\overline{\gamma}_{n}, &  j=-1,\\[1ex]
        0, & j=0,1,\ldots, n-1, \\[1ex]
        \gamma_{n}, & j=n
      \end{array} \right. \quad n \geq 1.
\end{equation*}
Here,   $\dsp \gamma_n = \frac{4 d_{n+1} }{(1+ic_{n+1})}\gamma_{n-1}$,\ $n\geq 1$,  \ with \  $\dsp \gamma_0 = \nu_0 = \frac{2d_1}{1+ic_1}$. \\[-2ex]

\end{lema}

\noindent {\bf Proof}.   Applying the respective three term recurrence formulas in %
\[
   U_{n}(z) = Q_{n}(z) R_{n-1}(z) - Q_{n-1}(z) R_{n}(z),  \quad n \geq 1,
\]
there follows \ $U_1(z) = 2d_1$  and
\begin{equation*} \label{Determinant-formula}
    U_{n+1}(z) = 4 d_{n+1} z U_{n}(z) = 2^{2n+1} d_1 d_2 \cdots d_{n+1} z^{n}, \quad n \geq 1.
\end{equation*}

Hence, considering the series expansions in terms of the origin and infinity,
\begin{equation} \label{Correspondence-1}
    \frac{Q_{n}(z)}{R_{n}(z)} - \frac{Q_{n-1}(z)}{R_{n-1}(z)} = \left\{
      \begin{array}{l}
        \dsp \frac{\overline{\gamma}_{n-1}}{\overline{r}_{n-1,n-1}}\, z^{n-1} + O\big(z^{n}\big), \\[3ex]
        \dsp \frac{\gamma_{n-1}}{r_{n-1,n-1}}\, \frac{1}{z^{n}} + O\big((1/z)^{n+1}\big),
      \end{array}
      \right.
      n \geq 1,
\end{equation}
where \
$ \dsp    \gamma_{n-1} = \frac{2^{2n-1} d_1 d_2 \cdots d_{n}}{r_{n,n}},  \ n \geq 0$ .  \ That is,  there exist formal series expansions $E_0$ and $E_{\infty}$, as in (\ref{Eq-Two-Series}),  such that  the correspondences in (\ref{origin-infinity-correspondence-1}) hold.   Since  $R_{n}$ and $Q_{n}$ are self inversive,  from the infinity correspondence we also have
\[
    z^{-1}\overline{E_{\infty}(1/\bar{z})} - \frac{Q_{n}(z)}{R_{n}(z)} = \frac{\overline{\gamma}_{n}}{\overline{r}_{n,n}}\, z^{n} + O\big(z^{n+1}\big),  \quad n \geq 0.
\]
Comparing this with the origin correspondence we then have the symmetry property $E_{0}(z) = z^{-1}\overline{E_{\infty}(1/\bar{z})} $, which is equivalent to (\ref{Eq-ModifiedMoments-Symmetry}).
The remaining result of the theorem follows by considering the systems of equations in  the coefficients of  $Q_n(z)= \sum_{j=0}^{n-1} q_{n,j}z^{j}$ and $R_{n}(z) = \sum_{j=0}^{n} r_{n,j}\,z^j$, which follow from (\ref{origin-infinity-correspondence-1}).  \hfill \qed

For other such contribution on three term recurrence formulas of the type (\ref{Eq-TTRR-Rn}) we also refer to \cite[Thm. 2.1]{Ismail-Masson-1995}.

For the next lemma  we restrict the sequence  $\{d_{n}\}_{n=1}^{\infty}$ in (\ref{Eq-TTRR-Rn}) to be such that $ \{d_{n+1}\}_{n=1}^{\infty} = \{d_{1,\,n}\}_{n=1}^{\infty} $ is a positive chain sequence. Clearly, Lemma \ref{Lemma-Orto-Rn} continues to hold as long as  $d_1 \neq 0$.  From  (\ref {Eq-TTRR-Rn}) we have $\frac{R_{n}(1)}{2R_{n-1}(1)} \big[1 - \frac{R_{n+1}(1)}{2R_{n}(1)}\big] = d_{1,\,n}$, $n \geq 1$.
Hence, $\{\mathfrak{m}_{n}\}_{n=0}^{\infty}$, with
\begin{equation}  \label{mParam-Rep}
    \mathfrak{m}_{n} = 1 - \frac{R_{{n+1}}(1)}{2R_{n}(1)}, \quad n \geq 0,
\end{equation}
is the minimal parameter sequence of the positive chain sequence $\{d_{1,n}\}_{n=1}^{\infty}$.

Note that in \cite{Castillo-Costa-Ranga-Veronese-2013} we have assumed the stronger restriction that $\{d_{n}\}_{n=1}^{\infty}$ is a positive chain sequence. It is important to note that if $\{d_{n}\}_{n=1}^{\infty}$ is a positive chain sequence then  so is $\{d_{1,\,n}\}_{n=1}^{\infty}$.  However, the inverse of the above affirmation is not always true. 

If $\{m_n\}_{n=0}^{\infty}$ is the minimal parameter sequence of $\{d_{n}\}_{n=1}^{\infty}$  then $\{m_{1,n}\}_{n=0}^{\infty}$, where $m_{1,n} = m_{n+1}$, $n \geq 0$,  is also a parameter sequence of $\{d_{1,\,n}\}_{n=1}^{\infty}$ and the minimal parameter sequence $\{\mathfrak{m}_n\}_{n=0}^{\infty}$  of  $\{d_{1,\,n}\}_{n=1}^{\infty}$ is such that $\mathfrak{m}_n < m_{1,n}$, $n \geq 0$.  Moreover,  if $\{M_n\}_{n=0}^{\infty}$ is the maximal parameter sequence of $\{d_{n}\}_{n=1}^{\infty}$  then $\{M_{1,n}\}_{n=0}^{\infty}$, where $M_{1,n} = M_{n+1}$, $n \geq 0$,  is exactly the maximal parameter sequence of $\{d_{1,\,n}\}_{n=1}^{\infty}$.  Hence, only when $\{d_{1,\,n}\}_{n=1}^{\infty}$ is a positive chain sequence with its maximal parameter sequence $\{M_{1,\,n}\}_{n=0}^{\infty}$  different from its minimal parameter sequence, i.e. $M_{1,0} > 0$, then the choice of $d_1$ such that $0 < d_1 \leq M_{1,0}$ makes $\{d_{n}\}_{n=1}^{\infty}$ also a positive chain sequence.

A positive chain sequence, for which the maximal parameter sequence is the same as the minimal parameter sequence, is said to determine its parameter uniquely.  We will  refer to such a chain sequence as a {\bf SPPCS} (Single Parameter Positive Chain Sequence).  By Wall's criteria (see \cite[p.\,101]{Chihara-book}) for maximal parameter sequence, the sequence $\{d_n\}_{n=1}^{\infty}$ is a SPPCS if and only if 
\[
     \sum_{n=1}^{\infty} \prod_{k=1}^{n} \frac{m_k}{1-m_k} \ = \infty,
\]
where $\{m_n\}_{n=0}^{\infty}$ is the minimal parameter sequence of $\{d_n\}_{n=1}^{\infty}$.

\begin{lema} [\cite{DimRan-2013}]\label{Lemma-SzegoKernel-InterlacingZeros}
Let the real sequences $\{c_{n}\}_{n=1}^{\infty}$ and $\{d_{n+1}\}_{n=1}^{\infty}$ be such that $\{d_{n+1}\}_{n=1}^{\infty}$ is also a positive chain sequence. Let $\{R_n\}$ be the sequence of polynomials obtained from   $(\ref{Eq-TTRR-Rn})$ with the use of $\{c_{n}\}_{n=1}^{\infty}$ and $\{d_{n+1}\}_{n=1}^{\infty}$.  Then the polynomial  $R_{n}(z)$ has  all its $n$ zeros simple and that these zeros lie on the unit circle $|z| = 1$. By denoting  the zeros of $R_{n}$ by $z_{n,j} = e^{i\theta_{n,j}}$, $j=1,2, \ldots,n$,  where $\theta_{n,j} < \theta_{n,j+1}$, then there holds the interlacing property
\[
    0 < \theta_{n+1,1} < \theta_{n,1} < \theta_{n+1,2} < \cdots < \theta_{n,n} < \theta_{n+1,n+1} < 2 \pi, \quad n \geq 1.
\]
Moreover, if \ $V_{n}(z) = R_{n}^{\prime}(z) R_{n-1}(z) - R_{n-1}^{\prime}(z) R_{n}(z)$, $n \geq 1$, then
\begin{equation*} \label{Wronskian-identity}
 z_{n,j}^{-(n-2)}(z_{n,j}-1)^{-1}\,V_{n}(z_{n,j})  > \,0, \quad j = 1, 2, \ldots, n, \quad n \geq 1.
\end{equation*} \\[-6ex]

\end{lema}

The initial part of this lemma was  established in \cite{DimRan-2013} with the use of the functions
\begin{equation} \label{Definition-Gn}
    G_{n}(x) = (4z)^{-n/2} R_{n}(z), \quad n \geq 0,
\end{equation}
given by the transformation $2x = z^{1/2}+z^{-1/2}$. This transformation, which  maps the points $z=e^{i\theta}$ onto the points $x = \cos(\theta/2)$, is referred to in \cite{Zhedanov-1998} as the DG transformation.

Clearly, the zeros of the function $G_{n}(x)$ in $[-1,1]$ are $x_{n,j} = \cos(\theta_{n,j}/2)$, $j=1,2, \ldots, n$.  The  associated Christoffel-Darboux functions or Wronskians
\begin{equation*} \label{Wronskian-for-General-Gn}
   W_{n}(x) = G_{n}^{\prime}(x)G_{n-1}(x)-G_{n-1}^{\prime}(x)G_{n}(x), \quad n \geq 1,
\end{equation*}
which not necessarily remain positive throughout $[-1,1]$, but satisfy at the zeros of $G_n(x)$ 
\begin{equation*} \label{Wronskian-for-General-Gn-Positiveness}
     W_{n}(x_{n,j}) >  0,  \quad j=1,2, \ldots, n \quad \mbox{and} \quad n \geq 1.
\end{equation*}
The last part of the lemma, proved in \cite{Castillo-Costa-Ranga-Veronese-2013}, can be obtained as follows.
From (\ref{Definition-Gn})
\[
     W_{n}(x) = \frac{(4z)^{-(n-1)}}{z-1}\left[2zV_{n}(z) -R_{n-1}(z)R_{n}(z)\right], \quad n \geq 1
\]
and, hence,
\begin{equation*} \label{Wronskian-identity}
 \frac{z_{n,j}^{-(n-2)}}{z_{n,j}-1}V_{n}(z_{n,j}) = 2^{2n-3}W_{n}(x_{n,j}) , \quad j = 1, 2, \ldots, n, \quad n \geq 1.
\end{equation*}

\subsection{Recovering the first positive measure}

Now we briefly state the results obtained in \cite{Castillo-Costa-Ranga-Veronese-2013} under the stronger restriction that $\{d_{n}\}_{n=1}^{\infty}$ is a positive chain sequence.
By considering the rational functions
\[
     \frac{\tilde{A}_{n}(z)}{\tilde{B}_{n}(z)}  = \frac{R_n(z) - Q_n(z)}{(z-1) R_n(z)}, \quad n \geq 1,
\]
there hold
\[
    \tilde{F}_{0}(z) - \frac{\tilde{A}_{n}(z)}{\tilde{B}_{n}(z)} = \dsp  \frac{\overline{\gamma}_{n}}{\overline{r}_{n,n}}\, z^{n} + O\big(z^{n+1}\big) \quad \mbox{and} \quad
     \tilde{F}_{\infty}(z) - \frac{\tilde{A}_{n}(z)}{\tilde{B}_{n}(z)} = \dsp - \frac{\gamma_{n}}{r_{n,n}}\, \frac{1}{z^{n+2}} + O\big((1/z)^{n+3}\big),
\]
for $n \geq 1$, where $\tilde{F}_{0}(z) = - \sum_{n=1}^{\infty} \tilde{\mu}_{n}\,z^{n-1}$ and $\tilde{F}_{\infty}(z) = \sum_{n=0}^{\infty} \tilde{\mu}_{-n}\,z^{-n-1}$,
with $\tilde{\mu}_{0} = 1$ and
\begin{equation*} \label{Eq-Moment-Relations}
      \tilde{\mu}_{n} = 1 + \sum_{j=1}^{n} \nu_{j},   \quad   \tilde{\mu}_{-n} = 1 - \sum_{j=1}^{n} \nu_{-j+1}, \quad n \geq 1.
\end{equation*}
Since $\nu_{n} = - \overline{\nu}_{-n+1}$, $n\geq 1$, one finds $\tilde{\mu}_{n} = \overline{\tilde{\mu}}_{-n}$, $n\geq 1$.

If one defines the moment functional $\mathcal{\tilde{M}}$ by $\mathcal{\tilde{M}}[z^{-n}] = \tilde{\mu}_{n}$, $n = 0, \pm1, \pm2, \ldots ,$ then
\begin{equation*} \label{Eq-Functional-Relations-1}
    \mathcal{\tilde{M}}[z^{-n}] =  1 - \mathcal{N}\big[\frac{\,1-z^{-n}}{1-z}\big], \quad n = 0, \pm1, \pm2, \ldots \ .
\end{equation*}
Reciprocally, $\mathcal{N}[z^{-n}] = \mathcal{\tilde{M}}[z^{-n}(1-z)]$, $n = 0, \pm1, \pm2, \ldots \ $.

For the rational functions $\tilde{A}_n/\tilde{B}_n$ there also hold the partial decompositions
\begin{equation} \label{Eq-Partial-DecompOld}
   \frac{\tilde{A}_n(z)}{\tilde{B}_n(z)} = \frac{\tilde{\lambda}_{n,0}}{z-1} + \sum_{j=1}^{n} \frac{\tilde{\lambda}_{n,j}}{z-z_{n,j}}, \quad n \geq 1,
\end{equation}
where $z_{n,j}$, $j=1,2, \ldots, n$, are the zeros of $R_n(z)$. The numbers  $\tilde{\lambda}_{n,j}$  are such that
\[
   \tilde{\lambda}_{n,j} > 0, \  j=0, 1, \ldots, n \quad \mbox{and} \quad \sum_{j=0}^{n} \tilde{\lambda}_{n,j} = 1,
\]
for $n \geq 1$. Precisely, $\tilde{\lambda}_{n,0} = 1 - \frac{Q_{n}(1)}{R_{n}(1)}$ is such that
\[
     1 - d_1 = \tilde{\lambda}_{1,0} > \tilde{\lambda}_{2,0} > \cdots > \tilde{\lambda}_{n,0} > \tilde{\lambda}_{n+1,0} > \cdots \qquad \mbox{and} \quad \lim_{n \to \infty} \tilde{\lambda}_{n,0} = M_0,
\]
where $\{M_n\}_{n=0}^{\infty}$ is the maximal parameter sequence of the positive chain sequence $\{d_n\}_{n=1}^{\infty}$, and
\[
   \tilde{\lambda}_{n,j} = \frac{U_{n}(z_{n,j})}{(1-z_{n,j})V_{n}(z_{n,j})} = \frac{2^{2n-1} d_1 d_2 \cdots d_{n}}{z_{n,j}^{-(n-2)}(z_{n,j}-1)^{-1}\,V_{n}(z_{n,j})} \frac{z_{n,j}}{(z_{n,j} - 1)(1 - z_{n,j})} > 0, \quad
\]
for $j =1, 2, \ldots n$ and $n \geq 1$.

Thus, with the application of the Helley selection theorem (see, for example, \cite{JoNjTh-1989})  to the Riemann-Stieltjes integral representation of (\ref{Eq-Partial-DecompOld}), we obtain that $\mathcal{\tilde{M}}(\zeta^{-n}) = \int_{\mathcal{C}}\zeta^{-n}d \tilde{\mu}(\zeta)$, $n=0,\pm1, \pm2, \ldots $.

Now by considering also some further results obtained in \cite{Castillo-Costa-Ranga-Veronese-2013} we can state the following theorem.

\begin{theo}  \label{Thm-Castillo-Costa-Ranga-Veronese-2013}
  Let the real sequences $\{c_{n}\}_{n=1}^{\infty}$ and $\{d_{n}\}_{n=1}^{\infty}$ be such that $\{d_{n}\}_{n=1}^{\infty}$ is also a positive chain sequence. Let $\{m_n\}_{n=0}^{\infty}$ and $\{M_n\}_{n=0}^{\infty}$ be the respective minimal and maximal parameter sequences of $\{d_{n}\}_{n=1}^{\infty}$.

  Associated with the sequences $\{c_{n}\}_{n=1}^{\infty}$ and $\{d_{n}\}_{n=1}^{\infty}$ there exists a nontrivial probability measure $\tilde{\mu}$ on the unit circle, with $M_0$ as its mass at $z=1$,  such that the  corresponding sequence of monic OPUC $\{\tilde{S}_n\}$  is given by
\[
   \tilde{S}_0(z) = 1 \quad \mbox{and} \quad \tilde{S}_n(z) \prod_{k=1}^{n}(1+ic_k) = R_{n}(z) - \,2(1-m_n)R_{n-1}(z), \quad n \geq 1.
\]
Here,  $R_n$ are the  polynomials obtained from   $(\ref{Eq-TTRR-Rn})$ with the use of $\{c_{n}\}_{n=1}^{\infty}$ and $\{d_{n+1}\}_{n=1}^{\infty}$.

In particular, for $\tilde{\mu} $ the Verblunsky coefficients $\tilde{\alpha}_{n-1} = -\overline{S_n(0)}$ are given by
\[
     \tilde{\alpha}_{n-1} =  \frac{1}{\tau_{n}} \,\frac{1-2m_n-ic_n}{1+ic_n},   \quad n \geq 1,
\]
where $\tau_{n} = \prod_{k=1}^{n} (1-ic_k)(1+ic_k)$, $n \geq 1$.

Moreover, if $\mathcal{N}$ is the moment functional associated with  $\{c_{n}\}_{n=1}^{\infty}$ and $\{d_{n}\}_{n=1}^{\infty}$, as described in Lemma \ref{Lemma-Orto-Rn}, then for any Laurent polynomial $\ell$,
\[
    \mathcal{N}[\ell(z)] = \int_{\mathcal{C}} \ell(\zeta) (1-\zeta)d \tilde{\mu}(\zeta) .
\]

\end{theo}

Observe that if the sequences $\{c_n\}_{n=1}^{\infty}$ and $\{d_n\}_{n=1}^{\infty}$  are taken to be the same as $\{c_n^{(2)}\}_{n=1}^{\infty}$ and $\{d_n^{(2)}\}_{n=1}^{\infty}$   then the measure $\tilde{\mu}(z)$ given by Theorem \ref{Thm-Castillo-Costa-Ranga-Veronese-2013} is the same as the measure $\mu(\delta;z)$ considered in the beginning of the section.

\setcounter{equation}{0}
\section{Recovering another positive measure } \label{Sec-Recovering-NewSzego}

With the sequence of polynomials $\{Q_{n}\}$ and $\{R_{n}\}$ given by (\ref{Eq-TTRR-Rn}) we now show that one can also obtain another nontrivial probability measure, say $\hat{\mu}$, under the conditions that the real sequences $\{c_{n}\}_{n=1}^{\infty}$ and $\{d_{n}\}_{n=1}^{\infty}$ are such that $d_1 \neq 0$ and $\{d_{n+1}\}_{n=1}^{\infty} = \{d_{1,\,n}\}_{n=1}^{\infty}$ is a positive chain sequence.

Let $\{\hat{A}_{n}\}$  be the sequence of polynomials  given by
\begin{equation} \label{Szego-RationalFunction}
     \hat{A}_{n}(z) = \frac{(1+ic_{1})}{2z}\big[R_n(z) + \frac{(1-ic_1)}{2d_1}(z-1)Q_n(z)\big],  \quad n \geq 1.
\end{equation}
With the observation made in  section \ref{Sec-RecentDevlop} that $R_n(0) = \prod_{k=1}^{n}(1-i c_k)$, $n \geq 1$,  $Q_1(0) = 2d_1$  and $Q_n(0) = 2d_1\prod_{k=2}^{n}(1-i c_k)$, $n \geq 2$, it is easily verified that $\hat{A}_n$ is a polynomial of degree $n-1$, with its leading coefficient equals to $\prod_{k=1}^{n}(1+ic_k)$.

Since
\[
   \frac{\hat{A}_{n+1}(z)}{R_{n+1}(z)} - \frac{\hat{A}_{n}(z)}{R_{n}(z)}= \frac{1+c_1^2}{4d_1}\, (1- z^{-1})\left[\frac{Q_{n+1}(z)}{R_{n+1}(z)} - \frac{Q_{n}(z)}{R_{n}(z)}  \right], \quad n \geq 1,
\]
from (\ref{Correspondence-1}) there hold
\begin{equation*} \label{Correspondence-2}
    \frac{\hat{A}_{n+1}(z)}{R_{n+1}(z)} - \frac{\hat{A}_{n}(z)}{R_{n}(z)} = \left\{
      \begin{array}{l}
        \dsp -\frac{(1+c_1^2)\overline{\gamma}_{n}}{4d_1\overline{r}_{n,n}}\, z^{n-1} + O\big(z^{n}\big), \\[3ex]
        \dsp \frac{(1+c_1^2)\gamma_{n}}{4d_1r_{n,n}}\, \frac{1}{z^{n+1}} + O\big((1/z)^{n+2}\big),
      \end{array}
      \right.
      n \geq 1.
\end{equation*}
Hence, there exist series expansions $\hat{F}_0(z)$ and $\hat{F}_{\infty}(z)$ such that
\begin{equation} \label{origin-correspondence-2}
    \hat{F}_{0}(z) - \frac{\hat{A}_{n}(z)}{R_{n}(z)} = -\frac{(1+c_1^2)\overline{\gamma}_{n}}{4d_1\overline{r}_{n,n}}\, z^{n-1} + O\big(z^{n}\big), \quad n \geq 1 \qquad  \quad
\end{equation}
and
\begin{equation} \label{infinity-correspondence-2}
    \hat{F}_{\infty}(z) - \frac{\hat{A}_{n}(z)}{R_{n}(z)} = \frac{(1+c_1^2)\gamma_{n}}{4d_1 r_{n,n}}\, \frac{1}{z^{n+1}} + O\big((1/z)^{n+2}\big), \quad n \geq 1.
\end{equation}
Setting
\[
  \hat{F}_{0}(z) = -\sum_{n=0}^{\infty} \hat{\mu}_{n+1} z^{n}    \quad \mbox{and}   \quad \hat{F}_{\infty}(z) = \sum_{n=1}^{\infty} \hat{\mu}_{-n+1} z^{-n},
\]
we then have from  (\ref{Szego-RationalFunction}), (\ref{origin-correspondence-2}) and (\ref{infinity-correspondence-2}), together with the results of Lemma \ref{Lemma-Orto-Rn}, that the numbers $\hat{\mu}_n$ satisfy
\begin{equation*} \label{Eq-Moment-Relations}
  \begin{array}{l}
      \dsp \hat{\mu}_{n} = \frac{1+c_1^2}{4d_1}(\nu_{n}-\nu_{n+1}),   \quad  \hat{\mu}_{-n} = \frac{1+c_1^2}{4d_1}(\nu_{-n}-\nu_{-n+1}), \quad n \geq 1,
  \end{array}
\end{equation*}
with $\hat{\mu}_0 = 1$. Since $\nu_{j} = - \overline{\nu}_{-j+1}$, $j \geq 1$, hence there hold
\[
  \hat{\mu}_{n} = \overline{\hat{\mu}}_{-n}, \quad n \geq 1.
\]
If we define the moment functional $\hat{\mathcal{M}}$ by $\mathcal{\hat{M}}[z^{-n}] = \hat{\mu}_{n}$, $n=0, \pm1, \pm2, \ldots$, then
\begin{equation} \label{Eq-Functional-Relations-2}
    \mathcal{\hat{M}}[z^{-n}] =  \frac{1+c_1^2}{4d_1} \mathcal{N}\big[z^{-n}-z^{-n-1}\big], \quad n=0, \pm1, \pm2, \ldots \ ,
\end{equation}
where $\mathcal{N}$ is the linear functional defined as in Lemma \ref{Lemma-Orto-Rn}.

Since $\hat{A}_n$ is of degree $n-1$, with the results of Lemma \ref{Lemma-SzegoKernel-InterlacingZeros} there holds the partial decomposition
\begin{equation*} \label{Eq-Partial-Decomp}
   \frac{\hat{A}_n(z)}{R_n(z)} = \sum_{j=1}^{n} \frac{\hat{\lambda}_{n,j}}{z-z_{n,j}} ,
\end{equation*}
where $z_{n,j} = e^{i\theta_{n,j}}$, $j=1,2, \ldots, n$, are the zeros of $R_n(z)$ and
\[
   \hat{\lambda}_{n,j} = \frac{\hat{A}_{n}(z_{n,j})}{R_{n}^{\prime}(z_{n,j})} = \frac{(1+c_1^2)(z_{n,j} - 1)}{4 d_1 z_{n,j}}\frac{Q_{n}(z_{n,j})}{R_{n}^{\prime}(z_{n,j})}, \quad j =1, 2, \ldots n.
\]
Since $Q_n(z_{n,j})/R_n^{\prime}(z_{n,j})$ can be written as $U_n(z_{n,j})/V_n(z_{n,j})$, again from Lemma \ref{Lemma-SzegoKernel-InterlacingZeros},
\[
     \hat{\lambda}_{n,j} =  \frac{2^{2n-3}(1+c_1^2)d_2 \ldots d_n}{z_{n,j}^{-(n-2)}(z_{n,j} - 1)^{-1}V_{n}(z_{n,j})} > 0, \quad j=1,2,\ldots, n.
\]
In addition to the positiveness of the elements $\hat{\lambda}_{n,j}$, $j=0,1, 2, \ldots, n$, by considering the limit of  $z\hat{A}_{n}(z)/R_{n}(z)$, as $z \to \infty$, we also have
\[
     \sum_{j=1}^{n} \hat{\lambda}_{n,j} = 1.
\]
Hence,  if the step functions $\hat{\psi}_{n}(e^{i\theta})$, $n \geq 1$, are defined on $[0, 2\pi]$ by
\[
        \hat{\psi}_n(e^{i\theta}) = \left\{
        \begin{array}{ll}
          0, & 0 < \theta \leq \theta_{n,1}, \\[1ex]
          \sum_{j=1}^{k} \hat{\lambda}_{n,j}, & \theta_{n,k} < \theta \leq \theta_{n,k+1},\quad k=1,2, \ldots, n-1, \\[1ex]
          1, & \theta_{n,n} < \theta \leq 2\pi.
        \end{array} \right.
\]
then from the definition of the Riemann-Stieltjes integrals
\begin{equation*} \label{Eq-RiemannStieltjes-Representation}
     \frac{\hat{A}_n(z)}{R_n(z)} = \int_{\mathcal{C}} \frac{1}{z - \zeta}\, d \hat{\psi}_{n}(\zeta), \quad n \geq 1.
\end{equation*}
Hence, by the application of the Helley selection theorem (see \cite{JoNjTh-1989}) there exists a subsequence $\{n_j\}$ such that $\hat{\psi}_{n_j}(e^{i\theta})$ converges to a bonded non-decreasing function, say $\hat{\mu}(e^{i\theta})$, in $[0, 2\pi]$.

From (\ref{origin-correspondence-2}) and (\ref{infinity-correspondence-2}), since
\[
   \int_{\mathcal{C}} d \hat{\psi}_{n}(\zeta) = 1 \quad \mbox{and} \quad \int_{\mathcal{C}} \zeta^{k} d \hat{\psi}_{n}(\zeta) = \hat{\mu}_{-k}, \ \ k =\pm1,\pm2, \ldots, \pm n-2,
\]
we also have that
\begin{equation} \label{NewMeasure-NewMoments}
   \int_{\mathcal{C}} d \hat{\mu}(\zeta) = 1 = \mathcal{\hat{M}}[1] \quad \mbox{and} \quad \int_{\mathcal{C}} \zeta^{n} d \hat{\mu}(\zeta) = \hat{\mu}_{-n}= \mathcal{\hat{M}}[z^n], \ \ n =\pm 1,\pm 2, \ldots\, .
\end{equation}

Now we can state the following theorem.

\begin{theo} \label{Thm-Main1}
  Given the real sequences $\{c_{n}\}_{n=1}^{\infty}$ and $\{d_{n+1}\}_{n=1}^{\infty}$, where $\{d_{n+1}\}_{n=1}^{\infty} = \{d_{1,\,n}\}_{n=1}^{\infty}$ is also a positive chain sequence then, associated with these, there exists a nontrivial probability measure $\hat{\mu}$ on the unit circle such that for any Laurent polynomial $\ell$,
\begin{equation} \label{Eq-New-Moment-Relation}
   \int_\mathcal{C} \ell(\zeta)d \hat{\mu}(\zeta) =  \frac{1+c_1^2}{4d_1} \mathcal{N}\big[\ell(z)(1-z^{-1})\big] .
\end{equation}
Here $\mathcal{N}$ is the moment functional associated with   $\{c_{n}\}_{n=1}^{\infty}$ and $\{d_{n}\}_{n=1}^{\infty}$, as described in Lemma $\ref{Lemma-Orto-Rn}$, where the additional term $d_1$ is chosen such that $d_1 \neq 0$.   The sequence of monic OPUC $\{\hat{S}_n\}$ with respect to $\hat{\mu}$  are
\begin{equation} \label{Eq-NewOPUC}
     \hat{S}_n(z) = \frac{R_{n+1}(z) - \,2(1-\mathfrak{m}_n)R_{n}(z)}{(z-1)\prod_{k=1}^{n+1}(1+ic_k)}, \quad n \geq 0,
\end{equation}
where $R_n$ are the polynomials obtained  from   $(\ref{Eq-TTRR-Rn})$ with the use of the two sequences $\{c_{n}\}_{n=1}^{\infty}$ and $\{d_{n+1}\}_{n=1}^{\infty}$.
In particular, the  Verblunsky coefficients $\hat{\alpha}_{n-1} = - \overline{\hat{S}_n(0)}$ with respect to $\hat{\mu}$  are
\[
      \hat{\alpha}_{n-1} = - \frac{1}{\tau_{n}}\,\frac{1-2\mathfrak{m}_n-ic_{n+1}}{1-ic_{n+1}}, \quad n \geq 1.
\]
Here, $\{\mathfrak{m}_n\}_{n=0}^{\infty}$, where $\mathfrak{m}_{n} = 1 - \frac{R_{{n+1}}(1)}{2R_{n}(1)}$, $n \geq 0$, is the minimal parameter sequence of $\{d_{1,\,n}\}_{n=1}^{\infty}$ and $\tau_{n} = \prod_{k=1}^{n} (1-ic_k)/(1+ic_k)$, $n \geq 1$.
\end{theo}

\noindent {\bf Proof}.  The required probability measure  $\hat{\mu}$  is the one established  in (\ref{NewMeasure-NewMoments}).  We then have from (\ref{Eq-Functional-Relations-2}) the moment functional relations in (\ref{Eq-New-Moment-Relation}).
From these, clearly,
\begin{equation} \label{Modified-Moment-Rep}
   \frac{1+c_1^2}{4d_1} \mathcal{N}\big[z^{-k}(1-z^{-1})P(z)\big] = \int_\mathcal{C}\zeta^{-k} P(\zeta) d \hat{\mu}(\zeta), \quad k=0, \pm1, \pm2, \ldots \ ,
\end{equation}
for any polynomial $P$.

Now, from  Lemma \ref{Lemma-Orto-Rn}, since
\begin{equation} \label{LOrthogonality-Rn-inSec4}
   \mathcal{N}\big[z^{-n+k} R_n(z)\big] = \left\{ \begin{array}{ll}
                                                                           \  0, &  k=0,1, \ldots, n-1, \\
                                                                           \  \gamma_n,  & k = n,
                                                                         \end{array} \right.
\end{equation}
where
\[
     \gamma_0 =  \frac{2d_1}{1+ic_1} \  \   \ \mbox{and} \  \  \   \gamma_n = \frac{4 d_{n+1} }{(1+ic_{n+1})}\gamma_{n-1},\  \ n\geq 1,
\]
by considering the monic  polynomial of degee $n$
\[
     \hat{S}_n(z) = \frac{1}{\prod_{k=1}^{n+1}(1 + ic_k)}\, \frac{1}{z-1} \big[R_{n+1}(z) - \frac{R_{n+1}(1)}{R_{n}(1)}R_{n}(z)\big],
\]
we obtain from (\ref{Modified-Moment-Rep}) that
\[
     \int_\mathcal{C}\zeta^{-k} \hat{S}_n(\zeta) d \hat{\mu}(\zeta) = \frac{1+c_1^2}{4d_1} \mathcal{N}\big[z^{-k}(1-z^{-1})\hat{S}_n(z)\big] =  \left\{
                                                                   \begin{array}{ll}
                                                                        0, &  k=0,1, \ldots, n-1, \\
                                                                        \hat{\kappa}_n^{-2},  & k = n,
                                                                    \end{array} \right.
\]
where $\hat{\kappa}_n^{-2} = (1-\mathfrak{m}_n)\prod_{k=1}^{n} \frac{4d_{k+1}}{1+c_{k+1}^2}$.
Hence, with the use of (\ref{mParam-Rep}), we obtain the results for the OPUC  and, consequently, also about the Verblunsky coefficients. \hfill \qed

Since $R_n^{\ast}(z) = R_n(z)$, $n \geq 0$, we also have from Theorem \ref{Thm-Main1} that
\begin{equation} \label{Eq-NewReciproOPUC}
    \hat{S}_n^{\ast}(z) = \frac{R_{n+1}(z) - \,2(1-\mathfrak{m}_n)zR_{n}(z)}{(1-z)\prod_{k=1}^{n+1}(1-ic_k)}, \quad n \geq 0.
\end{equation}

Now combining the results obtained in Theorem \ref{Thm-Castillo-Costa-Ranga-Veronese-2013} and the above theorem we can state the following.

\begin{theo}    \label{Thm-MainChain1}
 If $\{d_{1,\,n}\}_{n=1}^{\infty}$ is not a SPPCS  then the measure $\hat{\mu}$, derived as in Theorem \ref{Thm-Main1}, is such that the integral  $J = \int_{\mathcal{C}} \zeta (\zeta-1)^{-1} (1-\zeta)^{-1} d\hat{\mu}$ exists.

\end{theo}

\noindent {\bf Proof}.  First we recall that SPPCS refers to single parameter positive chain sequence. Now from Theorem \ref{Thm-Main1} the measure $\hat{\mu}$ is such that
\[
   \int_\mathcal{C} \ell(\zeta)d \hat{\mu}(\zeta) =  \frac{1+c_1^2}{4d_1} \mathcal{N}\big[\ell(z)(1-z^{-1})\big] ,
\]
where $\mathcal{N}$ is the moment functional associated with   $\{c_{n}\}_{n=1}^{\infty}$ and $\{d_{n}\}_{n=1}^{\infty}$, as described in Lemma $\ref{Lemma-Orto-Rn}$, where the additional term $d_1$ is chosen such that $d_1 \neq 0$.

Since the chain sequence $\{d_{1,\,n}\}_{n=1}^{\infty}$ does not determine its parameter uniquely we have $M_1 > 0$, where $\{M_{n+1}\}_{n=0}^{\infty} = \{M_{1,n}\}_{n=0}^{\infty}$ is the maximal parameter sequence of $\{d_{1,\,n}\}_{n=1}^{\infty}$.  Hence, if we chose $d_1 \leq  M_1$ then $\{d_{n}\}_{n=1}^{\infty}$ is also a positive chain sequence and from Theorem \ref{Thm-Castillo-Costa-Ranga-Veronese-2013}  there exists a probability measure $\tilde{\mu}$ on the unit circle such that
\[
    \mathcal{N}[\ell(z)] = \int_{\mathcal{C}} \ell(\zeta) (1-\zeta)d \tilde{\mu}(\zeta) .
\]
Hence,
\[
   \int_\mathcal{C} \ell(\zeta)d \hat{\mu}(\zeta) =  \frac{1+c_1^2}{4d_1} \int_\mathcal{C}  \ell(\zeta) (1-\zeta^{-1})(1-\zeta) d \tilde{\mu}(\zeta),
\]
from which we conclude the required result of the theorem.  \hfill \qed

If the integral $J = \int_{\mathcal{C}} \zeta (1-\zeta)^{-1}(\zeta-1)^{-1}d\hat{\mu}$  exists then so does the integral $\int_{\mathcal{C}} \zeta (\zeta-1)^{-1}d\hat{\mu}(\zeta)$.  Clearly, the inverse is not always true.

For the next theorem we assume that the measure $\hat{\mu}$ is such that the integral \linebreak $\int_{\mathcal{C}} \zeta (\zeta-1)^{-1}d\hat{\mu}(\zeta)$ exists or, better still,  that the principal value integral $I = \dashint_{\mathcal{C}} \zeta (\zeta-1)^{-1}d\hat{\mu}(\zeta)$ exists.

 Here, what we mean by $I = \dashint_{\mathcal{C}} \zeta (\zeta-1)^{-1}d\hat{\mu}(\zeta)$ is that
\[
       I = \lim_{\epsilon \to 0} \int_{\epsilon}^{2\pi - \epsilon} \frac{e^{i\theta}}{e^{i\theta}-1} d\hat{\mu}(e^{i \theta}) .
\]
Clearly, if $\int_{\mathcal{C}} \zeta (\zeta-1)^{-1}d\hat{\mu}(\zeta)$ exists then also the principal value integral $I = \int_{\mathcal{C}} \zeta (\zeta-1)^{-1}d\hat{\mu}(\zeta)$.  Since $\hat{\mu}$ is a probability measure it is easily verified that $I + \overline{I} = 2\mathcal{R}e(I) = 1$. Moreover, if $\hat{\mu}$ also satisfies the symmetry $d\hat{\mu}(e^{i(2\pi-\theta)}) = - d\hat{\mu}(e^{i\theta})$  then $I = 1/2$.

As an example, we consider the probability measure $\hat{\mu}$ given by
\[
      \int_{\mathcal{C}}\ell(\zeta) d \hat{\mu}(\zeta) = (1-\delta)\int_{\mathcal{C}} \ell(\zeta) \frac{1}{2\pi i\, \zeta} d \zeta  + \delta\, \ell(i),
\]
where $0 \leq \delta < 1$. Then
\[
   I = \dashint_{\mathcal{C}}\frac{\zeta}{\zeta -1 } d \hat{\mu}(\zeta) = (1 - \delta) \lim_{\epsilon \to 0} \frac{1}{2\pi}\int_{0+\epsilon}^{2\pi-\epsilon} \frac{e^{i\theta}}{e^{i\theta}-1} d \theta   \   + \ \delta\, \frac{i}{i-1}.
\]
This can be written as
\[
   I -  \ \delta\, \frac{i}{i-1} =  (1 - \delta)  \frac{1}{2\pi} \lim_{\epsilon \to 0} \Big[ \int_{0+\epsilon}^{\pi} \frac{e^{i\theta}}{e^{i\theta}-1} d \theta + \int_{0+\epsilon}^{\pi} \frac{e^{-i\theta}}{e^{-i\theta}-1} d \theta\Big] = (1 - \delta)  \frac{1}{2\pi} \lim_{\epsilon \to 0} \int_{0+\epsilon}^{\pi} d\theta,
\]
from which $I = 1/2 - i\delta/2$.

\begin{theo} \label{Thm-Main2}
  Let the associated  nontrivial probability measure $\hat{\mu}$ on the unit circle as given by Theorem  \ref{Thm-Main1} is such that  the principal value integral $I = \dashint_{\mathcal{C}} \zeta (\zeta-1)^{-1} d \hat{\mu}(\zeta)$ exists.   If  $\mathcal{N}$ is the moment functional associated with   $\{c_{n}\}_{n=1}^{\infty}$ and $\{d_{n}\}_{n=1}^{\infty}$, as described in Lemma $\ref{Lemma-Orto-Rn}$, where the additional term $d_1$ is chosen such that $d_1 \neq 0$, then for any Laurent polynomial $\ell$,
\[
   \frac{1+c_1^2}{4d_1} \mathcal{N}\big[\ell(z)\big] =  \dashint_\mathcal{C}\ell(\zeta) \frac{\zeta}{\zeta-1}d \hat{\mu}(\zeta) \  + \  i \,t \, \ell(1),
\]
where  $t = - \mathcal{I}m[(1+ic_1)I]$. In particular, if  $\{d_{1,n}\}_{n=1}^{\infty}$ is not a SPPCS then  $t = 0$.

\end{theo}

\noindent {\bf Proof}.   From (\ref{Eq-New-Moment-Relation}) and the symmetry property (\ref{Eq-ModifiedMoments-Symmetry}), clearly $\mathcal{N}$ is such that 
\[
   \frac{1+c_1^2}{4d_1} \mathcal{N}\big[\ell(z)\big] =  \dashint_\mathcal{C}\ell(z) \frac{\zeta}{\zeta-1}d \hat{\mu}(\zeta) \  + \  i \,t \, \ell(1),
\]
for all Laurent polynomials $\ell$, where $t$, $-\infty < t < \infty$.  Thus, form  (\ref{LOrthogonality-Rn-inSec4}), also
\begin{equation*} \label{Eq-Cauchy-Lorthogonality}
  \begin{array}{ll}
    \dsp \frac{1+c_1^2}{4d_1} \mathcal{N}\big[z^{-n+k}R_n(z)\big]  & \dsp \  =  \dashint_\mathcal{C} \zeta^{-n+k} R_{n}(\zeta) \frac{\zeta}{\zeta-1} d \hat{\mu}(\zeta)  \  +  \  i \,t R_{n}(1) \\[3ex]
     & \  =
            \left\{ \begin{array}{ll}
                 0, &  k=0,1,\ldots, n-1, \\[1ex]
                 \dsp \frac{1+c_1^2}{4d_1} \gamma_n, & k=n.
            \end{array} \right.
   \end{array}
\end{equation*}
We can find the value of $t$ from $\mathcal{N}\big[z^{-1}R_1(z)\big] = 0$.  Since  $\dashint_{\mathcal{C}}(\zeta-1)^{-1}d \hat{\mu}(\zeta) = - \overline{I} $, we then have
\[
      \mathcal{N}\big[z^{-1}R_1(z)\big]  = (1+ic_1)I - (1-ic_1)\overline{I} + 2i t = 0,
\]
which gives the required value for $t$.

Now if  $\{d_{1,n}\}_{n=1}^{\infty}$ is not a SPPCS, by  letting $d_1 \leq  M_1$, then from  $ \int_\mathcal{C} \ell(\zeta)d \hat{\mu}(\zeta) =  \frac{1+c_1^2}{4d_1} \int_\mathcal{C}  \ell(\zeta) \linebreak (1-\zeta^{-1})(1-\zeta) d \tilde{\mu}(\zeta)$ we have
\[
     I = \frac{1+c_1^2}{4d_1}  \int_{\mathcal{C}} (1-\zeta) d \tilde{\mu}(\zeta),
\]
 which gives $\mathcal{I}m(I) = \frac{1+c_1^2}{4d_1} \mathcal{I}m(\tilde{\mu}_1)$, where $\tilde{\mu}_1 = \int_{\mathcal{C}} \zeta^{-1} d \tilde{\mu}(\zeta) = \overline{\tilde{\mu}}_{-1}$.  On the other hand from $ \mathcal{N}[\ell(z)] = \int_{\mathcal{C}} \ell(\zeta) (1-\zeta)d \tilde{\mu}(\zeta)$ and $\mathcal{N}[\zeta^{-1}R_1(\zeta)] = 0$ we have
 $c_1 = -\frac{1+c_1^2}{2d_1} \mathcal{I}m(\tilde{\mu}_1)$.  Thus, recalling also that  $\mathcal{R}e(I) = 1/2$  we have $t = - \mathcal{I}m[(1+ic_1)I] = 0$.    \hfill \qed

From  (\ref{Eq-NewOPUC}) and (\ref{Eq-NewReciproOPUC}), by elimination of $R_{n+1}$,
\begin{equation} \label{Eq-NewPOPUC-1}
   \frac{1}{\prod_{k=1}^{n}(1+ic_k)}R_n(z) = \frac{1+ic_{n+1}}{2(1-\mathfrak{m}_n)} \, \Big[\hat{S}_n(z) + \tau_{n+1}\,\hat{S}_n^{\ast}(z)\Big], \quad n \geq 0.
\end{equation}
Similarly, by elimination of $R_{n}$,
\begin{equation} \label{Eq-NewPOPUC-2}
   \frac{1}{\prod_{k=1}^{n+1}(1+ic_k)} R_{n+1}(z) =   z\hat{S}_{n}(z) + \tau_{n+1}\,\hat{S}_{n}^{\ast}(z), \quad n \geq 0.
\end{equation}

Observe that the polynomials appearing on both sides of (\ref{Eq-NewPOPUC-1}) and (\ref{Eq-NewPOPUC-2}) are monic polynomials.  Hence, we can also write (\ref{Eq-NewPOPUC-1}) in the form
\begin{equation*} \label{Eq-NewPOPUC-11}
   \frac{1}{\prod_{k=1}^{n}(1+ic_k)}R_n(z) = \frac{1}{1 - \tau_{n+1}\hat{\alpha}_{n-1}} \, \Big[\hat{S}_n(z) + \tau_{n+1}\,\hat{S}_n^{\ast}(z)\Big], \quad n \geq 0.
\end{equation*}
Moreover, using the recurrence formula (\ref{Szego-A-RR}) for the OPUC $\hat{S}_n$ in the above expression
\begin{equation*} \label{Eq-NewPOPUC-21}
   \frac{1}{\prod_{k=1}^{n}(1+ic_k)}R_n(z) = \Big[z\hat{S}_{n-1}(z) + \frac{\tau_{n+1} -\overline{\hat{\alpha}}_{n-1}}{1-\tau_{n+1}\hat{\alpha}_{n-1}}\,\hat{S}_{n-1}^{\ast}(z)\Big], \quad n \geq 1.
\end{equation*}
Comparing this with (\ref{Eq-NewPOPUC-2}) we have
\(
     \tau_{n} = \frac{\tau_{n+1} -\overline{\hat{\alpha}}_{n-1}}{1-\tau_{n+1}\hat{\alpha}_{n-1}}, \  n \geq 1,
\)
or equivalently,
\[
     \tau_{n+1} = \frac{\tau_{n}  + \overline{\hat{\alpha}}_{n-1}}{1 + \tau_{n}\hat{\alpha}_{n-1}}, \quad n \geq 1.
\]
The above recurrence formula for $\tau_{n}$ is useful if we can find an expression for $\tau_{1}  $ in terms of the measure $\hat{\mu}$.

Assuming $I = \dashint_{\mathcal{C}}(\zeta-1)^{-1}\zeta d \hat{\mu}(\zeta)$ exists, from  $\mathcal{N}\big[z^{-1}R_1(z)\big] = 0$  and $\tau_1 = (1-ic_1)/(1+ic_1)$, the required expression for $\tau_1$ is
\[
   \tau_1 = \frac{I+it}{\overline{I}-it}.
\]

From $\tau_{n} = \prod_{k=1}^{n} (1-ic_k)/(1+ic_k)$, $n \geq 1$, and the recurrence formula for $\{\tau_n\}$,
\[
      c_1 = i \frac{\tau_1 - 1}{\tau_1 + 1}, \quad c_{n+1} = \frac{ \mathcal{I}m(\tau_{n+1}\hat{\alpha}_{n-1})} {1 - \mathcal{R}e(\tau_{n+1}\hat{\alpha}_{n-1})}
      = \frac{ \mathcal{I}m(\tau_{n}\hat{\alpha}_{n-1})} {1 + \mathcal{R}e(\tau_{n}\hat{\alpha}_{n-1})}, \  \ n \geq 1.
\]

From the expression for the Verblunsky coefficients in Theorem \ref{Thm-Main1},
\[
    1 + \tau_{n} \hat{\alpha}_{n-1} = \frac{2\mathfrak{m}_n}{1- ic_{n+1}} \quad \mbox{and}\quad 1 - \tau_{n+1} \hat{\alpha}_{n-1} = \frac{2(1-\mathfrak{m}_n)}{1+ ic_{n+1}} \quad n \geq 1.
\]
Hence,
\[
    \mathfrak{m}_{n} = \frac{1}{2}  \frac{|1 + \tau_{n}  \hat{\alpha}_{n-1}|^2}{1 + \mathcal{R}e(\tau_{n}  \hat{\alpha}_{n-1})} =  \frac{1}{2} \frac{1 - |\tau_{n+1}  \hat{\alpha}_{n-1}|^2}{1 - \mathcal{R}e( \tau_{n+1} \hat{\alpha}_{n-1})} , \quad n \geq 1,
\]
\[
   1 - \mathfrak{m}_{n} =  \frac{1}{2}  \frac{1 - |\tau_{n}  \hat{\alpha}_{n-1}|^2}{1 + \mathcal{R}e(\tau_{n} \hat{\alpha}_{n-1})} = \frac{1}{2}  \frac{|1 - \tau_{n+1}  \hat{\alpha}_{n-1}|^2}{1 - \mathcal{R}e(\tau_{n+1}  \hat{\alpha}_{n-1})} ,   \quad n \geq 1,
\]
and $d_{n+1} = (1-\mathfrak{m}_{n-1})\mathfrak{m}_{n}$, $n \geq 1$, with $\mathfrak{m}_{0} = 0$.
%

\setcounter{equation}{0}
\section{From OPUC to DG1POP } \label{Sec-OPUCtoDG1POP}

From what was presented  in the beginning of section  \ref{Sec-RecentDevlop}, given the sequence $\{S_n\}$ of OPUC  associated with some nontrivial probability measure $\mu$ on the unit circle then the sequence $\{zS_{n}(z) - \rho_{n}^{(2)}S_{n}^{\ast}(z)\}$ of para-orthogonal polynomials, where
\[
     \rho_{0}^{(2)}  = 1 \quad \mbox{and} \quad \rho_{n}^{(2)} = \frac{S_n(1)}{S_n^{\ast}(1)}
           = \frac{\rho_{n-1}^{(2)} - \overline{\alpha}_{n-1}}{1 - \rho_{n-1}^{(2)}   \alpha_{n-1}} ,  \   \ n \geq 1,
\]
satisfy the three term recurrence formula (\ref{Eq-TTRR-DG2Ex}).

By observing the results obtained in section \ref{Sec-Recovering-NewSzego} we are now able to state the following theorem which gives  a general result about how to generate other sequences of para-orthogonal polynomials that satisfy three term recurrence formulas.

\begin{theo} \label{Thm-Main3}
Let $\{\alpha_n\}_{n=0}^{\infty}$ be such that $|\alpha_n| < 1$, $n\geq 0$, and let $\mu$ and $\{S_n\}$ be the respective nontrivial probability measure and sequence of OPUC for which $\alpha_n$ are the Verblunsky coefficients.  With $\rho_0$ chosen such that $|\rho_{0}| = 1$,  $\rho_{0} \neq 1$, let the sequence of numbers $\{\rho_n\}_{n=0}^{\infty}$  be given by
\begin{equation} \label{Eq-RhoRecurence1}
    \rho_{n} = \frac{\rho_{n-1} - \overline{\alpha}_{n-1}}{1 - \rho_{n-1} \alpha_{n-1}}, \quad n \geq 1,
\end{equation}
where  $\alpha_{n-1} = - \overline{S_n(0)}$, $n \geq 0$, are the Verblunsky coefficients.  Then the sequence  $\{R_n(z)\}$ of POPUC  given by
\[
   R_n(z)  \prod_{k=0}^{n-1} \frac{ 1 + \mathcal{R}e(\rho_{k}\alpha_{k-1})}{1 + \overline{\rho}_{k} \overline{\alpha}_{k-1}}  = zS_{n-1}(z) - \rho_{n-1}S_{n-1}^{\ast}(z), \quad n \geq 1,
\]
satisfy the three term recurrence formula
\begin{equation} \label{Eq-TTRR-DG1Ex}
   R_{n+1}(z) = \big[(1+ic_{n+1})z + (1-ic_{n+1})\big] R_{n}(z) - 4\,d_{n+1} z R_{n-1}(z), \quad n \geq 1,
\end{equation}
with $R_{0}(z) = 1$ and $R_{1}(z) = (1+ic_{1})z + (1-ic_{1})$. Here,  the real sequence $\{c_n\}$ and the positive chain sequence $\{d_{n+1}\}_{n=1}^{\infty} = \{d_{1,n}\}_{n=1}^{\infty} $ are such that
\begin{equation*}\label{Eq-CoeffsTTRR-2}
     c_{1} = i\frac{\rho_0 + 1}{\rho_0 - 1}, \quad c_{n+1} =  \frac{ -\mathcal{I}m(\rho_{n}\alpha_{n-1})} {1 + \mathcal{R}e(\rho_{n}\alpha_{n-1})}   \quad \mbox{and}  \quad d_{1,n} = (1-\mathfrak{m}_{n-1})\mathfrak{m}_{n}, \quad n \geq 1,
\end{equation*}
where  $\{\mathfrak{m}_n\}_{n=0}^{\infty}$ is the minimal parameter sequence  of $\{d_{1,n}\}_{n=1}^{\infty}$ given by
\[
      \mathfrak{m}_{n} = \frac{1}{2} \frac{1 - |\rho_{n}  \alpha_{n-1}|^2}{1 + \mathcal{R}e( \rho_{n} \alpha_{n-1})}, \quad n \geq 0.
\]
The measure $\mu$ is such that the integral $J = \int_{\mathcal{C}} \zeta (1-\zeta)^{-1}(\zeta-1)^{-1}d\mu$  exists if and only if  there exists a $\rho_0$ $(|\rho_0| = 1, \rho_0 \neq 1)$ such that the corresponding positive chain sequence $\{d_{1,n}\}_{n=1}^{\infty}$ is not a SPPCS.

\end{theo}

\noindent {\bf Proof}.  We consider the monic polynomials $P_0(z) = 1$ and $P_{n} (z) = zS_{n-1}(z) - \rho_{n-1}S_{n-1}^{\ast}(z)$, $n \geq 1$, and show that, as in \cite[Thm. 2.1]{Costa-Felix-Ranga-2013},
\begin{equation} \label{Eq-MonicRecurrence}
    P_{n+1} (z) = (z  + \frac{\rho_{n}}{\rho_{n-1}})  P_{n} (z)   - (1 - \overline{\rho}_{n-1}\,\overline{ \alpha}_{n-1})(1 + \rho_{n-1} \alpha_{n-2}) z P_{n-1} (z), \quad n\geq 1.
\end{equation}
With the use of  (\ref{Szego-A-RR})  and (\ref{Eq-RhoRecurence1}), first observe that
\[
    P_{n} (z) = \frac{S_n(z) - \rho_n S_n^{\ast}(z)}{1 + \rho_n \alpha_{n-1}}, \quad n \geq 0.
\]
Hence for the polynomials
\[
   p_{n+1}(z) = P_{n+1}(z) + (1 - \overline{\rho}_{n-1}\,\overline{ \alpha}_{n-1})(1 + \rho_{n-1} \alpha_{n-2}) z P_{n-1} (z), \quad n\geq 1,
\]
we have
\[
    p_{n+1}(z) = [zS_n(z) - \rho_n S_n^{\ast}(z)] + (1 - \overline{\rho}_{n-1}\,\overline{ \alpha}_{n-1}) z [S_{n-1}(z) - \rho_{n-1}S_{n-1}^{\ast}(z)], \quad n\geq 1.
\]
Once again using (\ref{Szego-A-RR})  and then (\ref{Eq-RhoRecurence1}) this reduces to
\[
    p_{n+1}(z) = z [zS_{n-1}(z) - \rho_{n-1}S_{n-1}^{\ast}(z)] + \frac{\rho_{n}}{\rho_{n-1}} [zS_{n-1}(z) - \rho_{n-1}S_{n-1}^{\ast}(z)] ,  \quad n \geq 1,
\]
thus proving (\ref{Eq-MonicRecurrence}).  From (\ref{Eq-MonicRecurrence}), with the observation
\[
      (1 + \rho_{n} \alpha_{n-1}) (1 - \rho_{n-1}\alpha_{n-1}) =  (1 + \overline{\rho}_{n} \overline{\alpha}_{n-1}) (1 - \overline{\rho}_{n-1}\,\overline{ \alpha}_{n-1}) = 1 - |\alpha_{n-1}|^2, \quad n \geq 0,  
\]
one can easily derive (\ref{Eq-TTRR-DG1Ex}).

Since, $\mathfrak{m}_{0} = 0$, $\mathfrak{m}_{n} > 0$, $n\geq 1$, and
\[
        1 -\mathfrak{m}_{n} = \frac{1}{2} \frac{|1 + \rho_{n}  \alpha_{n-1}|^2}{1 + \mathcal{R}e( \rho_{n} \alpha_{n-1})} > 0, \quad n \geq 1,
\]
one can also conclude that $\{d_{1,n}\}_{n=1}^{\infty}$ is a positive chain sequence and that $\{\mathfrak{m}_n\}_{n=0}^{\infty}$ is  its minimal parameter sequence.

Before proving the last part of the theorem, observe that 
\[ 
   \rho_{n} = -\prod_{k=1}^{n+1}(1-ic_k)/(1+ic_k) = -\tau_{n+1}, \quad n \geq 0,  
\]    
where $\tau_{n}$ are the same as in Theorem \ref{Thm-Main1}.     Hence, starting from the three term recurrence formula (\ref{Eq-TTRR-DG1Ex}) we find that the measure $\hat{\mu}$  given by Theorem \ref{Thm-Main1}  is the same as the measure $\mu$.   Therefore, if $\{d_{1,n}\}_{n=1}^{\infty}$ is not a SPPCS  then, as in Theorem \ref{Thm-MainChain1}, the measure $\mu$ is such that the integral $J$ exists.

On the other hand, if $\mu$ is such that $J$ exists we can consider the probability measure $\tilde{\mu}$ given by
\[
       d \tilde{\mu} (\zeta) = J^{-1} \frac{\zeta}{(\zeta-1)(1-\zeta)} d \mu(\zeta).
\]
Hence, for the para-orthogonal polynomials $(z-1)R_n^{(2)}( z) = (z-1)R_n^{(2)}(\tilde{\mu}, z)$ we obtain from the results presented in section \ref{Sec-RecentDevlop} that
\[
   R_{n+1}^{(2)}(z) = \big[(1+ic_{n+1}^{(2)})z + (1-ic_{n+1}^{(2)})\big] R_{n}^{(2)}(z) - 4\,d_{n+1}^{(2)} z R_{n-1}^{(2)}(z),
\]
with $R_{0}^{(2)}(z) = 1$ and $R_{1}^{(2)}(z) = (1+ic_{1}^{(2)})z + (1-ic_{1}^{(2)})$, where $\{c_n^{(2)}\}_{n=1}^{\infty}$ is a real sequence and $\{d_{n+1}^{(2)}\}_{n=1}^{\infty}$ is a positive chain sequence.  Clearly,   $\{d_{n+1}^{(2)}\}_{n=1}^{\infty}$ is also not a SPPCS. Consequently, using results obtained in Theorems \ref{Thm-Main1} and \ref{Thm-MainChain1}, we conclude that with the choice $\rho_0$  such that $c_1^{(2)} = i (\rho_0 + 1)/(\rho_0-1)$  the resulting positive chain sequence $\{d_{1,n}\}_{n=1}^{\infty}$ = $\{d_{n+1}^{(2)}\}_{n=1}^{\infty}$  is not a SPPCS.    \hfill \qed

Note that every choice of $\rho_0$ generates a different sequence of polynomials $\{R_n\}$ which we can consider as a different extension of the DG1POP. By observing also the results establish in Theorem \ref{Thm-Main2}, we can  state the following Theorem which gives the orthogonality property associated  with each of these  the extension to the DG1POP.

\begin{theo} \label{Thm-Main4}
  Let $\mu$ be a nontrivial probability measure on the unit circle such that  the principal value integral   $I = \dashint_\mathcal{C} \zeta(\zeta-1)^{-1} d \mu(\zeta)$  exists and  let  $\{S_n\}$ be the  sequence of monic OPUC with respect to $\mu$.  Let $t$ be such that $-\infty < t < \infty$ and let the sequence  $\{\rho_{n}(t)\}_{n=0}^{\infty}$, with the property $|\rho_{n}(t)| = 1$, be given by
\[
	     \dsp \rho_{0}(t)=  - \frac{I+it}{\overline{I} -it}\quad \mbox{and} \quad  \rho_{n}(t) = \frac{\rho_{n-1}(t) - \overline{\alpha}_{n-1}}{1 - \rho_{n-1}(t) \alpha_{n-1}}, \quad n \geq 1,
\]
where $\alpha_{n-1} = - \overline{S_n(0)} $, $n \geq 0$.  Then for the sequence  $\{R_n(t;z)\}$ of POPUC, in $z$,  given by
\[
   \quad R_n(t;z)  \prod_{k=0}^{n-1} \frac{ 1 + \mathcal{R}e(\rho_{k}(t)\alpha_{k-1})}{1 + \overline{\rho_{k}(t)\alpha_{k-1}}}  = zS_{n-1}(z) - \rho_{n-1}(t)S_{n-1}^{\ast}(z), \quad n \geq 1,
\]
the following three term recurrence formula hold.
\begin{equation*}
   R_{n+1}(t;z) = \big[(1+ic_{n+1}(t))z + (1-ic_{n+1}(t))\big] R_{n}(t;z) - 4\,d_{n+1}(t) z R_{n-1}(t;z),
\end{equation*}
with $R_{0}(t;z) = 1$ and $R_{1}(t;z) = (1+ic_{1}(t))z + (1-ic_{1}(t))$,  where the real sequences $\{c_n(t)\}$ and $\{d_{n+1}(t)\}$ are such that
\begin{equation*}\label{Eq-CoeffsTTRR-2}
  \begin{array}l
     \dsp c_{1}(t) = i\frac{\rho_0(t) + 1}{\rho_0(t) - 1}   =  -2\big(t + \mathcal{I}m(I)\big)  \quad \mbox{and} \\ [2ex]
     \dsp c_{n+1}(t)=  \frac{ -\mathcal{I}m\big(\rho_{n}(t)\alpha_{n-1}\big)} {1 + \mathcal{R}e\big(\rho_{n}(t)\alpha_{n-1}\big)}, \ \  d_{n+1}(t) = \big(1-\mathfrak{m}_{n-1}(t)\big)\mathfrak{m}_{n}(t), \ \ n \geq 1,
   \end{array}
\end{equation*}
where   $\{\mathfrak{m}_n(t)\}_{n=0}^{\infty}$ is the minimal parameter sequence of the positive chain sequence $\{d_{n+1} (t)\}_{n=1}^{\infty}$,  which is given by
\[
      \mathfrak{m}_{n}(t) = \frac{1}{2} \frac{1 - |\rho_{n}(t)  \alpha_{n-1}|^2}{1 + \mathcal{R}e\big( \rho_{n}(t) \alpha_{n-1}\big)}, \quad n \geq 0.
\]
The sequence $\{d_{n+1} (t)\}_{n=1}^{\infty}$ is always a SPPCS if $t \neq 0$.  However,  $\{d_{n+1} (0)\}_{n=1}^{\infty}$ is not a SPPCS if and only if the integral $J = \int_\mathcal{C} \zeta(\zeta-1)^{-1} (1-\zeta)^{-1} d \mu(\zeta)$  exists.

With respect to the measure $\mu$, the polynomials $R_n(t;z)$ satisfy the L-orthogonality \linebreak $\mathcal{N}^{(t)}[z^{-n+k}R_n(t;z)] = 0$, $k=0,1,\ldots, n-1$,
with respect to the moment fuctional $\mathcal{N}^{(t)}$ given by
\begin{equation} \label{Eq-TheoMain4}
     \frac{1 + 4(t+\mathcal{I}m(I))^2}{4d_1} \, \mathcal{N}^{(t)}[\ell(z)]  = \dashint_{\mathcal{C}} \ell(\zeta) \frac{\zeta}{\zeta-1} d\mu (\zeta) + i\,t\, \ell(1).
\end{equation}
\end{theo}

Observe that $I$ in the above theorem is such that $I + \overline{I} = \mu_0 = 1$.  Hence, $\mathcal{R}e(I) = 1/2$ and $\rho_{0}(t)  \neq 1$ if $-\infty < t < \infty$.

A further observation is that the L-orthogonality $\mathcal{N}^{(t)}[z^{-n+k}R_n(t;z)] = 0$, $k=0,1,\ldots, n-1$,  holds even without the constant multiple $(1 + 4(t+\mathcal{I}m(I))^2)/(4d_1)$ in the definition (\ref{Eq-TheoMain4}) of $\mathcal{N}^{(t)}$.  However, if we maintain this constant multiple then $\mathcal{N}^{(t)}[1]  = 2d_1/(1+ic_1(t))$  and the modified moments $\nu_n = \mathcal{N}^{(t)}[z^{-n}]$ , $n =0, \pm 1, \pm 2, \ldots$, will be exactly the same as those follow from Lemma \ref{Lemma-Orto-Rn}.

When the measure  $\mu$ satisfy the  symmetric $d\mu(e^{i\theta}) = - d\mu(e^{i (2\pi - \theta)})$ then the Verblunsky coefficients $\alpha_{n-1}$ are all real and the principal value integral $I = \dashint_{\mathcal{C}} \zeta(\zeta - 1)^{-1} d \mu(\zeta)$ always exists and takes the value $1/2$.  Consequently, $\rho_{n}(0) = -1$, $n \geq 1$, and using  Wall's criteria  for maximal parameter sequences (see the formula prior to Lemma \ref{Lemma-SzegoKernel-InterlacingZeros}) we can state the following.

\begin{coro}
    Given the  sequence of real numbers $\{\alpha_{n}\}_{n=0}^{\infty}$, such that $|\alpha_{n}| < 1$, $n \geq 0$, let $\mu$ be the positive measure on the unit circle for which  $\alpha_{n}$ are its Verblunsky coefficients.  Then $\mu$ is such that integral $J = \int_\mathcal{C} \zeta(\zeta-1)^{-1} (1-\zeta)^{-1} d \mu(\zeta)$ exists if and only if the infinite series
 \[
       \sum_{n=1}^{\infty}  \prod_{k=1}^n\frac{1 + \alpha_{k-1}}{1 - \alpha_{k-1}},
 \]
 is convergent.
\end{coro}

\setcounter{equation}{0}
\section{Some further applications} \label{Sec-FurtherConseq}

First we see how one can obtain information about the OPUC associated with the positive measure $\hat{\mu}$, given by  $d \hat{\mu}(\zeta) = \zeta^{-1}(1-\zeta)(\zeta-1)d\mu$, if one knows information about the orthogonal polynomials associated with the measure $\mu$.

\begin{theo}  \label{Thm-Applic1}
    Let $\{S_n\}$ be the sequence of monic OPUC with respect to the nontrivial   measure $\mu$.  Let
\begin{equation} \label{Eq1-Thm-Applic1}
   c_{n} =  \frac{-\mathcal{I}m (\rho_{n-1}^{(2)}\alpha_{n-1})} {1-\mathcal{R}e(\rho_{n-1}^{(2)}\alpha_{n-1})}
   \quad \mbox{and} \quad
   g_{1, n-1} = \frac{1}{2} \frac{\big|1 - \rho_{n-1}^{(2)} \alpha_{n-1}\big|^2}{\big[1 - \mathcal{R}e(\rho_{n-1}^{(2)}\alpha_{n-1})\big]},  \quad n \geq 1,
\end{equation}
where   \   $\rho_0^{(2)}  = 1$,
\[ \dsp
    \alpha_{n-1} = - \overline{S_n(0)}  \quad  \mbox{and}  \quad   \rho_n^{(2)}
           = \frac{\rho_{n-1}^{(2)} - \overline{\alpha}_{n-1}}{1 - \rho_{n-1}^{(2)}   \alpha_{n-1}} , \quad   n \geq 1.  \
\]
Then the Verblunsky coefficients $\hat{\alpha}_{n-1} = - \overline{\hat{S}_n(0)}$  associated with the measure $\hat{\mu}$, where
\[
     d \hat{\mu}(\zeta) =  \  \frac{(\zeta-1)(1-\zeta)}{\zeta}d\mu(\zeta),
\]
are given by
\begin{equation}  \label{Eq2-Thm-Applic1}
        \hat{\alpha}_{n-1} = - \frac{1}{\rho_{n}^{(2)}}\,\frac{1-2\mathfrak{m}_n - ic_{n+1}}{1-ic_{n+1}}, \quad n \geq 1,
\end{equation}
where  \
\( \dsp
    \mathfrak{m}_{0} = 0 \  \mbox{and}  \   \mathfrak{m}_{n} = \frac{(1-g_{1,n-1})}{(1-\mathfrak{m}_{n-1})}\,g_{1,n},   \   n \geq 1.
\)

\end{theo}

\noindent {\bf Proof}.  With the  $\{c_n\}_{n=1}^{\infty}$ and $\{d_n\}_{n=1}^{\infty}$, where $d_1 = g_{1,0}$ and $d_{n+1} = (1-g_{1,n-1}) g_{1,n}$, $n \geq 1$, given by (\ref{Eq1-Thm-Applic1}), let $\mathcal{N}$  be the moment functional given by Lemma \ref{Lemma-Orto-Rn}.  Then  from results given in section \ref{Sec-RecentDevlop} that 
\[
     \mathcal{N}[\ell(z)] = a \int_{\mathcal{C}} \ell(\zeta) (1-\zeta) d \mu(\zeta),
\]
for any Laurent polynomial $\ell$.   Here, $a$ is some nonzero constant.  Observe that in section \ref{Sec-RecentDevlop}, $\mathcal{N}$ is normalized such that if $\mu$ is a probability measure than $a = 1$. 

On the other hand from results given in section  \ref{Sec-Recovering-NewSzego} we have the following. Associated with $\{c_n\}_{n=1}^{\infty}$ and $\{d_{n+1}\}_{n=1}^{\infty}$ there exists a nontrivial probability measure $\hat{\mu}$ on the unit circle such that 
\[
      \int_{\mathcal{C}} \ell(\zeta)  d \hat{\mu}(\zeta) = \frac{1+c_1^2}{4d_1}  \mathcal{N}[\ell(z)(1-z^{-1})]  ,
\]
for any Laurent polynomial $\ell$.    Moreover, the Verblunsky coefficients associated with the measure $\hat{\mu}$ are exactly as given in (\ref{Eq2-Thm-Applic1}).

From the two representations for $\mathcal{N}$ we have
\[
      \int_{\mathcal{C}} \ell(\zeta)  d \hat{\mu}(\zeta) = a \frac{1+c_1^2}{4d_1} \int_{\mathcal{C}} \ell(\zeta) (1-\zeta^{-1})(1-\zeta) d \mu(\zeta).
\]
Thus, establishing the proof of the theorem.    \hfill \qed

Thus, the above theorem is obtained by first using the results obtained in section \ref{Sec-RecentDevlop} and then followed by using the results given in section \ref{Sec-Recovering-NewSzego}.  Now, if we start by using the results given in sections \ref{Sec-Recovering-NewSzego} and \ref{Sec-OPUCtoDG1POP} and then  use the results given in section \ref{Sec-RecentDevlop} we obtain the following  theorem.

\begin{theo}  \label{Thm-Applic2}
    Let $\mu$ be a positive measure on the unit circle such that the integral \linebreak $J = \int_{\mathcal{C}} \zeta (1-\zeta)^{-1}(\zeta-1)^{-1}d\mu$ exists and let  $\{S_n\}$  be the associated sequence of monic OPUC.  Let the sequence of numbers $\{\rho_{n}\}_{n=0}^{\infty}$ be such that
\[
     \dsp \rho_{0} =  -\frac{I}{\,\overline{I}\,} \quad \mbox{and} \quad  \rho_{n}
               = \frac{\rho_{n-1}  - \overline{\alpha_{n-1}}}{1 - \rho_{n-1}\alpha_{n-1}}, \quad n \geq 1,
\]
where $I = \int_\mathcal{C} (\zeta-1)^{-1} d \mu(\zeta)$ and  $\alpha_{n-1} = - \overline{S_n(0)} $, $n \geq 0$.   Let the real  sequences $\{c_n\}_{n=1}^{\infty}$ and $\{d_{1,n}\}_{n=1}^{\infty}$ be such that
\begin{equation*}\label{Eq-CoeffsTTRR-2}
     c_{1} = i\frac{\rho_0 + 1}{\rho_0 - 1} \quad \mbox{and} \quad  c_{n+1}=  \frac{ -\mathcal{I}m(\rho_{n}\alpha_{n-1})} {1 + \mathcal{R}e(\rho_{n}\alpha_{n-1})} , \  \  d_{1,n} = (1 - \mathfrak{m}_{n-1}) \mathfrak{m}_{n},  \quad n\geq 1,
\end{equation*}
where $\{\mathfrak{m}_{n}\}_{n=0} ^{\infty}$ is the minimal parameter sequence of the chain sequence $\{d_{n+1}\}_{n=1}^{\infty}$ given by
\[
    \mathfrak{m}_{n} = \frac{1}{2} \frac{1 - |\rho_{n}  \alpha_{n-1}|^2}{1 + \mathcal{R}e( \rho_{n} \alpha_{n-1})}, \quad n \geq 0.
\]
Let $\{M_{n}\}_{n=1}^{\infty}$ be the maximal parameter sequence of $\{d_{n+1}\}_{n=1}^{\infty}$.  The existence of the integral $J$ assures that $M_{1} > 0$.

With $0 \leq  M_0 < 1$, let
\[
    m_1 = (1-M_0)M_{1} \quad \mbox{and} \quad m_{n+1} = d_{n+1}/m_{n}, \  \  n \geq 1.
\]
Then
\[
     \tilde{\alpha}_{n-1} =  - \frac{1}{\rho_{n-1}} \,\frac{1-2m_n-ic_n}{1+ic_n},   \quad n \geq 1,
\]
are the Verblunsky coefficients associated with the probability measure $\tilde{\mu}$ given by
\[
     \int_{\mathcal{C}} f(\zeta) d \tilde{\mu} = M_0 f(1) +  \frac{1-M_0}{J} \int_{\mathcal{C}} f(\zeta) \frac{\zeta}{(1-\zeta)(\zeta-1)} d\mu(\zeta).
\]

\end{theo}

\setcounter{equation}{0}
\section{Examples} \label{Sec-Examples}

\noindent {\bf Example 1.}  As a first example we consider the real sequence $\{c_n\}_{n=1}^{\infty}$ and the positive chain sequence $\{d_{n+1}\}_{n=1}^{\infty}$, where
\begin{equation} \label{Eq-Example1}
      c_1 = c, \quad c_{n+1} = 0,  \   n \geq 1, \qquad d_2 = \frac{1}{2}, \quad d_{n+2} = \frac{1}{4}, \  n  \geq 1,
\end{equation}
with $c$ real.   From (\ref{Eq-TTRR-Rn}) one can easily obtain that
\[
     R_n(z) = (1+ic)z^{n} + (1-ic) \quad \mbox{and} \quad Q_n(z) = 2d_1 \frac{z^{n} - 1}{z-1}, \quad n \geq 1,
\]
where $d_1 \neq 0$ is arbitrary. Thus, we can verify  that
\[
   \begin{array}{l}
      \dsp -\sum_{k=0}^{\infty} \frac{-2d_1}{1-ic} z^k - \frac{Q_n(z)}{R_n(z)} = \frac{\overline{\gamma}_nz^n}{\overline{r}_{n,n}} + O(z^{n+1}) ,  \\[3ex]
      \dsp  \ \  \, \sum_{k=1}^{\infty} \frac{2d_1}{1+ic} z^{-k} - \frac{Q_n(z)}{R_n(z)} = \frac{\gamma_n}{r_{n,n}z^{n+1}} + O((1/z)^{n+2}) ,
   \end{array}  \ \ \ n \geq 0,
\]
with $\dsp \gamma_0 = \frac{2d_1}{1+ic}$, \  \  $\dsp \gamma_n = \frac{4d_1}{1+ic}$, \ $n \geq 1$, \   \  $r_{0,0} = 1$ \  and \ $r_{n,n} = 1 + ic$, \ $n \geq 1$.

As defined in  Lemma \ref{Lemma-Orto-Rn}, thus the moment functional $\mathcal{N}$ associated with the sequences $\{c_n\}_{n=1}^{\infty}$ and $\{d_{n}\}_{n=1}^{\infty}$  is given by
\[
      \mathcal{N}[z^{-n}] = \nu_n =  \frac{-2d_1}{1-ic}, \quad \mathcal{N}[z^{n-1}] = \nu_{-n+1} = \frac{2d_1}{1+ic}, \quad n \geq 1.
\]

Now for the polynomials $\hat{A}_n$ defined in section \ref{Sec-Recovering-NewSzego} we have
\begin{equation*} \label{Ex-Szego-RationalFunction}
     \hat{A}_{n}(z) = \frac{(1+ic_{1})}{2z}\big[R_n(z) + \frac{(1-ic_1)}{2d_1}(z-1)Q_n(z)\big] = (1+ic)z^{n-1},  \quad n \geq 1.
\end{equation*}
Consequently, with the series expansions $\hat{F}_{0}(z) = -\sum_{n=0}^{\infty} \hat{\mu}_{n+1} z^{n}$ and $\hat{F}_{\infty}(z) = \sum_{n=1}^{\infty} \hat{\mu}_{-n+1} z^{-n}$, where $\hat{\mu}_0 = 1$ and $\hat{\mu}_n = 0$, $n \geq 1$, we can verify that

\[
    \hat{F}_{0}(z) - \frac{\hat{A}_{n}(z)}{R_{n}(z)} = -\frac{(1+c^2)\overline{\gamma}_{n}}{4d_1\overline{r}_{n,n}}\, z^{n-1} + O\big(z^{n}\big), \quad n \geq 1 \qquad  \quad
\]
and
\[
    \hat{F}_{\infty}(z) - \frac{\hat{A}_{n}(z)}{R_{n}(z)} = \frac{(1+c^2)\gamma_{n}}{4d_1 r_{n,n}}\, \frac{1}{z^{n+1}} + O\big((1/z)^{n+2}\big), \quad n \geq 1.
\]

Since the zeros $z_{n,j}$ of $R_n$ are the $n^{th}$ roots of $(1-ic)/(1+ic)$,  we also directly obtain that
\[
        \frac{\hat{A}_{n}(z)}{R_{n}(z)} = \sum_{j=1}^{n} \frac{1/n}{z - z_{n,j}}, \quad n \geq 1.
\]
Observing that the zeros $z_{n,j}$, $j=1,2,\ldots, n$,  are of equal distance from each other in the unit circle, we  conclude that the probability measure $\hat{\mu}$ given by Theorem \ref{Thm-Main1} is $d\hat{\mu}(e^{i\theta}) = (2\pi)^{-1}d\theta$.   Note that we can also write $d\hat{\mu}(\zeta) = (2\pi i\zeta)^{-1}d\zeta$.

The measure $d\hat{\mu}(\zeta) = (2\pi i\zeta)^{-1}d\zeta$ is well known  (in Simon \cite{Simon-book-p1} referred to as the free case) and the associated monic OPUC are $S_n(z) = z^n$, $n
\geq 0$.

Since the measure $d\hat{\mu}(\zeta) = (2\pi i\zeta)^{-1}d\zeta$  is also such that $d\hat{\mu}(e^{i(2\pi-\theta)}) = - d\hat{\mu}(e^{i\theta})$,  we have $I = \linebreak \dashint_{\mathcal{C}}(\zeta-1)^{-1}\zeta d \hat{\mu}(\zeta) = 1/2$.  Hence, from Theorem \ref{Thm-Main4}
\[
     \rho_{n}(t) = - \frac{1 +i 2t}{1-i2t} , \quad n \geq 0
\]
and the polynomials $\{R_n(t;z)\}$ given by
\[
   \quad R_n(t;z)  \prod_{k=0}^{n-1} \frac{ 1 + \mathcal{R}e(\rho_{k}(t)\alpha_{k-1})}{1+\overline{\rho_{k}(t) \alpha_{k-1}}}  = zS_{n-1}(z) - \rho_{n-1}(t)S_{n-1}^{\ast}(z), \quad n \geq 1,
\]
satisfy the three term recurrence formula
\begin{equation*}
   R_{n+1}(t;z) = \big[(1+ic_{n+1}(t))z + (1-ic_{n+1}(t))\big] R_{n}(t;z) - 4\,d_{n+1}(t) z R_{n-1}(t;z),
\end{equation*}
with $R_{0}(t;z) = 1$ and $R_{1}(t;z) = (1+ic_{1}(t))z + (1-ic_{1}(t))$,  where the real sequences $\{c_n(t)\}$ and $\{d_{n+1}(t)\}$ are such that
\[
      c_1(t) = -2t, \quad c_{n+1}(t)= 0,  \   n \geq 1, \qquad d_2(t) = \frac{1}{2}, \quad d_{n+2}(t) = \frac{1}{4}, \  n  \geq 1.
\]
The L-orthogonality property $ \tilde{\mathcal{N}}^{(t)}\big[z^{-n+k} R_{n}(t;z)\big] = 0$, $0 \leq k \leq n-1$, holds with the moment functional $\tilde{\mathcal{N}}^{(t)}$ given by
\[
   \tilde{\mathcal{N}}^{(t)}\big[\ell(z)\big] =  \dashint_\mathcal{C}\ell(\zeta) \frac{\zeta}{\zeta-1} \, \frac{1}{2\pi i \zeta}d \zeta \  + \  i \,t\, \ell(1),
\]
for any Laurent polynomial $\ell$.

Hence, the polynomials $R_n$ that follow from (\ref{Eq-Example1}) satisfy the L-orthogonality
\[
    \tilde{\mathcal{N}}^{(-c/2)}\big[z^{-n+k} R_{n}(z)\big] = 0,  \quad 0 \leq k \leq n-1
\]
and
\[
    \nu_n = \mathcal{N}\big[z^{-n}\big] = \frac{4d_1}{1+c^2}\tilde{\mathcal{N}}^{(-c/2)}\big[z^{-n}\big], \quad n = 0, \pm1, \pm2, \ldots \ .
\]

Since the integral $J = \int_{\mathcal{C}} \zeta (1-\zeta)^{-1}(\zeta-1)^{-1}d\hat{\mu}$ does not exist,  the  sequence $\{d_{n+1}(t)\}_{n=0}^{\infty}$ is a SPPCS for all $t$, as confirmed by Theorem \ref{Thm-Main4}.

\quad \vspace{0ex}

\noindent {\bf Example 2.}  Here we start with the probability measure $\mu$ given by
\[
     \int_{\mathcal{C}} \ell(\zeta) d \mu(\zeta) =  \int_{\mathcal{C}}  \ell(\zeta) \frac{(1-\zeta)(\zeta -1)}{4\pi i \zeta^2} d \zeta ,
\]
and derive the polynomials $R_n(t;z)$ and their L-orthogonality property given by Theorem \ref{Thm-Main4}.  Again, many information regarding the OPUC associated with the measure $\mu$ are well known.

For example, the associated Verblunsky coefficients are
\[
      \alpha_{n-1} = -\overline{S_n(0)} = - \frac{1}{n+1}, \quad n \geq 1.
\]
Clearly the integral $J = \int_{\mathcal{C}}(1-\zeta)^{-1} (\zeta - 1)^{-1} \zeta \,d \mu(\zeta)$ exists and that it is easy to show that
\[
       I =  \dashint_{\mathcal{C}} \frac{\zeta}{\zeta - 1} d \mu(\zeta)  = \int_{\mathcal{C}} \frac{\zeta}{\zeta - 1} d \mu(\zeta)  = \frac{1}{2}.
\]
We also verify by induction that the sequence $\{\rho_{n}(t)\}$ generated by
\[
	     \dsp \rho_{0}(t)=  -\frac{I+it}{\overline{I} -it}\quad \mbox{and} \quad  \rho_{n}(t) = \frac{\rho_{n-1}(t)  - \overline{\alpha_{n-1}}}{1 -  \rho_{n-1}(t) \alpha_{n-1}}, \quad n \geq 1,
\]
can be explicitly given as
\[
   \rho_{n}(t)  = -\frac{1 + i (n+1)(n+2)t}{1 - i (n+1)(n+2)t}, \quad n\geq 0.
\]
With this we obtain
\[
     c_{n+1}(t)= - \frac{2(n+1)t}{1 + n(n+1)^3 t^2} \quad   \mbox{and}  \quad \mathfrak{m}_{n}(t) = \frac{n}{2(n+1)} \frac{1 + (n+1)^2 (n+2)^2 t^2}{1 + n (n+1)^2 (n+2)t^2}, \quad n \geq 0.
\]

For the sequence of polynomials $\{R_{n}(t;z)\}$ given by
\begin{equation*}
   R_{n+1}(t;z) = \big[(1+ic_{n+1}(t))z + (1-ic_{n+1}(t))\big] R_{n}(t;z) - 4\,d_{n+1}(t) z R_{n-1}(t;z),
\end{equation*}
with $R_{0}(t;z) = 1$ and $R_{1}(t;z) = (1+ic_{1}(t))z + (1-ic_{1}(t))$,  the L-orthogonality  property $\tilde{\mathcal{N}}^{(t)}[z^{-n+k}R_n(t;z)] = 0$, $k=0,1,\ldots, n-1$, holds
with respect to the moment functional
\[
      \tilde{\mathcal{N}}^{(t)}[\ell(z)]  =  \int_{\mathcal{C}} \ell(\zeta) \frac{1 - \zeta}{4 \pi i \,\zeta} d\zeta + i\,t\, \ell(1).
\]

Since,
\[
      1 - \mathfrak{m}_{n}(t) = \frac{n+2}{2(n+1)} \frac{1 + n^2 (n+1)^2 t^2}{1 + n (n+1)^2 (n+2)t^2},  \quad n \geq 1,
\]
we have
\[
     \sum_{n =1}^{\infty}  \prod_{k=1}^{n} \frac{\mathfrak{m}_{k}(t)}{1 - \mathfrak{m}_{k}(t)} =  \frac{2}{1 + 4t^2}\sum_{n =1}^{\infty}  \Big[\frac{1}{(n+1)(n+2)} + (n+1)(n+2) t^2 \Big].
\]
The above infinite series converge if $t = 0$ and diverge otherwise.  Hence, using Wall's criteria (see Chihara \cite[p. 101]{Chihara-book}), the positive chain sequence $\{d_{n+1}(t)\} = \{d_{1,n}(t)\}$, where
\[
      d_{1,n}(t) = \big(1-\mathfrak{m}_{n-1}(t)\big)\mathfrak{m}_{n}(t), \quad n \geq 1,
\]
is not a SPPCS  if $t=0$ and is a SPPCS  otherwise.  Again, confirming the results in Theorem \ref{Thm-Main4}.

\quad \vspace{0ex}

\noindent {\bf Example 3.}  Let the real sequences $\{c_n\}_{n=1}^{\infty}$ and $\{d_{n+1}\}_{n=1}^{\infty}$ be  given by
\begin{equation*} \label{Eq-Special-Example}
 \begin{array}l
  \dsp c_n = \frac{\eta}{\lambda+n}, \quad d_{n+1} = \frac{1}{4} \frac{n(2\lambda+n+1)}{(\lambda+n)(\lambda+n+1)}, \ \ n\geq 1,
 \end{array}
\end{equation*}
where \ $\lambda, \eta \in \mathbb{R}$  and  $\lambda > -1$.

Observe that $d_{n+1} = d_{1,n} = (1-\mathfrak{m}_{n-1})\mathfrak{m}_{n}$, $n \geq 1$,  where
\begin{equation} \label{Eq-Special-ParamSeq}
   \quad \mathfrak{m}_{n} = \frac{n}{2(\lambda+n+1)}, \quad \ n \geq 0.
\end{equation}
Hence,  with our assumption $\lambda > -1$, the sequence $\{d_{1,n}\}_{n=1}^{\infty}$ is a positive chain sequence with $\{\mathfrak{m}_n \}_{n=0}^{\infty}$ as its minimal parameter sequence.

The polynomials $R_n$ obtained from the above sequences $\{c_n\}$ and $\{d_{n+1}\}$, together with the recurrence formula (\ref{Eq-TTRR-Rn}), are
\begin{equation} \label{Eq-HypergeometricRn}
      R_n(z) = \frac{(2\lambda+2)_n}{(\lambda+1)_n} \, _2F_1(-n,b+1;\,b+\overline{b}+2;\,1-z),  \quad n \geq 1,
\end{equation}
where  $b = \lambda + i \eta$.

With $d_1 \neq 0$ arbitrary, by considering the series expansion of $\,_2F_1(1,-b;\,\overline{b}+2;\,z)$ via its Gauss type continued fraction  expansion, we have
\[
    \frac{Q_{n}(z)}{R_{n}(z)} \sim \frac{2d_1}{1-ic_1}\,_2F_1(1,-b;\,\overline{b}+2;\,z)  = -\nu_{1} - \nu_{2} z - \nu_{3}z^2 - \ldots \, ,
\]
from which
\begin{equation} \label{Eq-Special-Modified-Moments}
     \nu_n = d_1 \frac{b + \overline{b}+2}{b+1} \frac{(-b-1)_{n}}{(\overline{b}+1)_n}, \quad n \geq 1.
\end{equation}
Since $\nu_{n} = - \overline{\nu}_{-n+1}$, with the convention $(a)_n = \Gamma(a+n)/\Gamma(a)$ for all integer values of $n$, the   above expression for $\nu_n$ is also valid for $n \leq 0$.

For the moment functional $\mathcal{N}$ defined by  $\mathcal{N}[\zeta^{-n}]  = \nu_n$, $n=0,\pm1,\pm2, \ldots$, we have
\begin{equation} \label{Eq-Integral-Representation-for-SMM}
  \begin{array}{l}
    \dsp  \mathcal{N}[\ell(z)]  =
      \int_{\mathcal{C}}\ell(\zeta) \,d \psi(b;\,\zeta),  \\[2ex]
  \end{array}
\end{equation}
for all Laurent polynomials $\ell$, where \
\(
  \dsp d \psi(b;\,z) = d_1 \frac{i\,|\Gamma(b+1)|^2 (b+ \overline{b} + 2)}{2\pi\, \Gamma(b+\overline{b}+2)} (-z)^{-\overline{b}-1} (1-z)^{b+\overline{b}+1} dz .
\)

\noindent Here, the branch cuts in $(-z)^{-\overline{b}} = (e^{-i\pi}z)^{-\overline{b}}$ and $(1-z)^{b+\overline{b}} = \big(e^{-i\pi}(z-1)\big)^{b+\overline{b}}$ are along the positive real axis.

For the proofs of the above results, starting from (\ref{Eq-HypergeometricRn}), we cite for example \cite{Castillo-Costa-Ranga-Veronese-2013} and \cite{Costa-Felix-Ranga-2013}.  However, considering what was necessary,  it was assumed in  \cite{Castillo-Costa-Ranga-Veronese-2013} and \cite{Costa-Felix-Ranga-2013} that $\lambda > -1/2$. We remark that when  $\lambda > -1/2$ the maximal parameter sequence $\{M_{1,n}\}_{n=0}^{\infty}$ of $\{d_{1,n}\}_{n=1}^{\infty}$ is
\[
      M_{1,n} = \frac{1}{2}\frac{2\lambda +n +1}{\lambda+n+1},  \quad n \geq 0.
\]
Now, the nontrivial probability measure $\hat{\mu}$ that follows from  Theorem \ref{Thm-Main1} is
\begin{equation} \label{Eq-NewMeasureExample3}
     d \hat{\mu}(\zeta) =  \frac{|\Gamma(b+2)|^2 }{2\pi\, \Gamma(b+\overline{b}+3)} (-z)^{-\overline{b}-2} (1-z)^{b+\overline{b}+2} dz.
\end{equation}
Observe that (\ref{Eq-NewMeasureExample3}) can also be written in the equivalent form
\[
    d \hat{\mu}(e^{i\theta})= \frac{2^{b+\overline{b}+2}|\Gamma(b+2)|^2 }{2\pi\, \Gamma(b+\overline{b}+3)} e^{(\pi-\theta)\mathcal{I}m(b+1)} [\sin^{2}(\theta/2)]^{\mathcal{R}e(b+1)} d\theta . \qquad
\]
The moments associated with $\hat{\mu}$ are
\[
      \overline{\mu}_{-n} = \mu_{n}  = \int_{\mathcal{C}} \zeta^{-n} d \hat{\mu}(\zeta) = \frac{1+c_1^2}{4d_1} [\nu_n - \nu_{n+1}] = \frac{(-b-1)_n}{(\overline{b}+2)_n}, \quad n = 0, 1, 2, \ldots \ .
\]
Hence, from results established in \cite{Ranga-2010} the associated monic OPUC and Verblunsky coefficients are
\[
    \begin{array}l 
          \dsp \hat{S}_n(z) = \frac{(b+\overline{b}+3)_n}{(b+2)_n} \, _2F_1(-n,b+2;\,b+\overline{b}+3;\,1-z),  \\[3ex]
          \dsp \hat{\alpha}_{n-1} = - \frac{(b+1)_{n}}{(\overline{b}+2)_{n}},
    \end{array}   n \geq 1.
\]
The above expression for $\hat{S}_n$ can also be easily derived from (\ref{Eq-NewOPUC}) and (\ref{Eq-HypergeometricRn}).

It is important to clarify that the above hypergeometric polynomials $R_n$ and $\hat{S}_n$ are polynomials belonging to a family of biorthogonal polynomials mentioned in the comments of Prof. Richard Askey on the Gabor Szeg\H{o}: Collected Papers \cite[p.\,304]{Askey-Szego}.  However, that $\hat{S}_n$ are OPUC with respect to the nontrivial probability  measure in (\ref{Eq-NewMeasureExample3}) and also that $R_n$ are para-orthogonal polynomials were first observed in \cite{Ranga-2010}.  See also the comments given at the end of \cite{Ranga-2010}.

Recall that we have assumed $\lambda > -1$ and  for the integral $I = \int_{\mathcal{C}}(\zeta-1)^{-1}\zeta\, d \hat{\mu}(\zeta)$  we have from  (\ref{Eq-Special-Modified-Moments}), (\ref{Eq-Integral-Representation-for-SMM})  and (\ref{Eq-NewMeasureExample3}) that
\[
       I = \frac{\overline{b}+1}{b+\overline{b}+2}.
\]
Let $t$ be such that $-\infty < t < \infty$ and let the sequence  $\{\rho_{n}(t)\}_{n=0}^{\infty}$, with the property $|\rho_{n}(t)| = 1$, be given by
\begin{equation} \label{Eq-RhoExample3}
	     \dsp \rho_{0}(t)=  -\frac{I+it}{\overline{I} -it}\quad \mbox{and} \quad  \rho_{n}(t) = \frac{\rho_{n-1}(t)  - \overline{\alpha}_{n-1}}{1 - \rho_{n-1}(t) \alpha_{n-1}}, \quad n \geq 1,
\end{equation}
where $\alpha_{n-1} = \hat{\alpha}_{n-1} = - \overline{\hat{S}_n(0)} $, $n \geq 0$.  Then from Theorem \ref{Thm-Main4} for the sequence  $\{R_n(t;z)\}$ of POPUC, in $z$,  given by
\[
   \quad R_n(t;z)  \prod_{k=0}^{n-1} \frac{ 1 + \mathcal{R}e(\rho_{k}(t)\alpha_{k-1})}{1 +  \overline{\rho_{k}(t) \alpha_{k-1}}}  = z\hat{S}_{n-1}(z) - \rho_{n-1}(t)\hat{S}_{n-1}^{\ast}(z), \quad n \geq 1,
\]
the following three term recurrence formula hold.
\begin{equation*}
   R_{n+1}(t;z) = \big[(1+ic_{n+1}(t))z + (1-ic_{n+1}(t))\big] R_{n}(t;z) - 4\,d_{n+1}(t) z R_{n-1}(t;z),
\end{equation*}
with $R_{0}(t;z) = 1$ and $R_{1}(t;z) = (1+ic_{1}(t))z + (1-ic_{1}(t))$,  where the real sequences $\{c_n(t)\}$ and $\{d_{n+1}(t)\}$ are such that
\begin{equation*}\label{Eq-CoeffsTTRR-2}
      c_{1}(t) = i\frac{\rho_0(t) + 1}{\rho_0(t) - 1} \quad \mbox{and} \ \  c_{n+1}(t)=  \frac{ \mathcal{I}m(\rho_{n}(t)\alpha_{n-1})} {1 - \mathcal{R}e(\rho_{n}(t)\alpha_{n-1})},  \ \  \  d_{n+1}(t) = \big(1-\mathfrak{m}_{n-1}(t)\big)\mathfrak{m}_{n}(t), \quad n \geq 1,
\end{equation*}
where   $\{\mathfrak{m}_n(t)\}_{n=0}^{\infty}$ is the minimal parameter sequence of the the chain sequence $\{d_{n+1} (t)\}_{n=1}^{\infty}$,  which is given by
\[
      \mathfrak{m}_{n}(t) = \frac{1}{2} \frac{1 - |\rho_{n}(t)  \alpha_{n-1}|^2}{1 + \mathcal{R}e( \rho_{n}(t) \alpha_{n-1})}, \quad n \geq 0.
\]
With respect to the measure $\mu$, the polynomials $R_n(t;z)$ satisfy the L-orthogonality \linebreak $\mathcal{N}^{(t)}[z^{-n+k}R_n(t;z)] = 0$, $k=0,1,\ldots, n-1$,
with respect to the moment fuctional $\mathcal{N}^{(t)}$ given by
\begin{equation*}
     \frac{1 + 4(t+\mathcal{I}m(I))^2}{4d_1} \, \mathcal{N}^{(t)}[\ell(z)]  = \dashint_{\mathcal{C}} \ell(\zeta) \frac{\zeta}{\zeta-1} d\mu (\zeta) + i\,t\, \ell(1).
\end{equation*}

Clearly, the sequence of polynomials $\{R_n(0;z)\}_{n=0}^{\infty}$ and the chain sequence $\{d_{n+1}(0)\}_{n=1}^{\infty}$ obtained above  are respectively the same as the sequences $\{R_n(z)\}_{n=0}^{\infty}$ and  $\{d_{n+1}\}_{n=1}^{\infty}$.  In fact, from (\ref{Eq-RhoExample3}) we can also easily verify that $\rho_{n}(0)  = - (\overline{b}+1)_{n+1}/ (b+1)_{n+1}$, $n \geq 0$. 

Not that when $\lambda > -1/2$ the maximal parameter sequence $\{M_{1,n}\}_{n=0}^{\infty}$ of $\{d_{n+1}(0)\}_{n=1}^{\infty}$ is different from its minimal parameter sequence $\{\mathfrak{m}_n\}_{n=0}^{\infty}$.   Since the integral $J = \int_{\mathcal{C}} \zeta (1-\zeta)^{-1}(\zeta-1)^{-1}d\hat{\mu}$ exists when $\lambda > -1/2$, this is exactly what we expect from Theorem \ref{Thm-Main4}.  However, when $-1/2 \geq \lambda > -1$, the integral $J = \int_{\mathcal{C}} \zeta (1-\zeta)^{-1}(\zeta-1)^{-1}d\hat{\mu}$ does not exist and, hence from Theorem \ref{Thm-Main4}, the minimal parameter sequence of $\{d_{n+1}(0)\}_{n=1}^{\infty}$ is also its maximal parameter sequence. This we can verify as follows.

From (\ref{Eq-Special-ParamSeq}) we have
\[
    \sum_{n=1}^{\infty} \Big[ \prod_{k=1}^{n}\frac{\mathfrak{m}_k(0)}{1 - \mathfrak{m}_k(0)}\Big] =  \sum_{n=1}^{\infty} \Big[ \prod_{k=1}^{n}\frac{\mathfrak{m}_k}{1 - \mathfrak{m}_k}\Big]  = \sum_{n=1}^{\infty} \frac{(1)_n}{(2\lambda+3)_n}.
\]
The infinite series on the right hand side  can also be written as $_2F_1(1,1;\,2\lambda+3;\,1)$, which is (see \cite[p. 62]{Andrews-Book}) convergent for $\lambda > -1/2$ and divergent for $-1/2 \geq \lambda > -1$.

\end{document}